\documentclass[12pt]{article}
\usepackage{epsfig,wrapfig,epic}
\usepackage{amsmath,latexsym}
\usepackage{amscd}
\usepackage{mathtools}
\usepackage{rotating}
\usepackage{amsfonts}
\usepackage{mathrsfs}
\usepackage{stmaryrd,bbm}
\usepackage{shadow}
\usepackage{lineno}
\usepackage{setspace}
\usepackage{color}
\newcommand{\blue}[1]{\textcolor{blue}{#1}}

\newcommand{\mapform}[5]{\begin{array}{ccrcl}
#1 & : & #2 & \longrightarrow & #3 \\
& & #4 & \longmapsto & #5
\end{array}}

\usepackage{dsfont}
\usepackage{graphics}
\usepackage{graphicx}
\usepackage{amssymb}

\title{A Newton's Iteration Converges Quadratically to Nonisolated 
Solutions Too}

\author{Zhonggang Zeng \thanks{Department of Mathematics,
Northeastern Illinois University, Chicago, Illinois 60625, USA.
~~email:~{\tt zzeng@neiu.edu}. 
~Research is supported in part by NSF under grant DMS-1620337.}}
\DeclareMathSymbol{\bdC}{\mathbin}{AMSb}{'103}
\DeclareMathAlphabet{\mathpzc}{OT1}{pzc}{m}{it}
\DeclareMathSymbol{\mP}{\mathbin}{AMSb}{'120}

\begin{document}
\newcommand{\rkr}{{\mbox{\scriptsize rank-$r$}}}
\newcommand{\rk}[1]{{\mbox{\scriptsize rank-{#1}}}}

\newtheorem{theorem}{Theorem}[section]
\newtheorem{lemma}[theorem]{Lemma}

\newtheorem{definition}[theorem]{Definition}
\newtheorem{example}[theorem]{Example}

\newtheorem{remark}[theorem]{Remark}
\newtheorem{corollary}[theorem]{Corollary}

\numberwithin{equation}{section}

\newtheorem{algrthm}{Algorithm}
\newtheorem{Problem}{PROBLEM}

\newcommand{\bdx}{\mathbf{x}}
\newcommand{\bdy}{\mathbf{y}}
\newcommand{\bdz}{\mathbf{z}}
\newcommand{\bdf}{\mathbf{f}}
\newcommand{\bdg}{\mathbf{g}}
\newcommand{\bdo}{\mathbf{0}}
\newcommand{\bde}{\mathbf{e}}
\newcommand{\F}{\mathbbm{F}}
\newcommand{\R}{\mathbbm{R}}
\newcommand{\C}{\mathbbm{C}}
\newcommand{\bdb}{\mathbf{b}}
\newcommand{\cK}{\mathpzc{Kernel}}
\newcommand{\cR}{\mathpzc{Range}}
\newcommand{\cM}{{\mathcal M}}
\newcommand{\rank}[1]{\mathpzc{rank}\left(\,#1\,\right)}
\newcommand{\ranka}[2]{\mathpzc{rank}_{#1}\left(\,#2\,\right)}
\newcommand{\h}{{{\mbox{\tiny $\mathsf{H}$}}}}
\newcommand{\sg}{\sigma}
\newcommand{\bdu}{\mathbf{u}}
\newcommand{\bdv}{\mathbf{v}}
\newcommand{\al}{\alpha}
\newcommand{\bt}{\beta}
\newcommand{\la}{\lambda}
\newcommand{\qed}{${~} $ \hfill \raisebox{-0.3ex}{\LARGE $\Box$}}
\newcommand{\Dl}{\Delta}
\newcommand{\dl}{\delta}
\newcommand{\dm}{\mathpzc{dim}}
\newcommand{\nullity}[1]{\mathpzc{nullity}\left(\,#1\,\right)}
\newcommand{\eps}{\varepsilon}
\newcommand{\spn}{\mathpzc{span}}


\date{}

\maketitle

\centerline{\em \small
This work is dedicated to the memory of Tien-Yien Li (1945-2020)}

\begin{abstract}
The textbook Newton's iteration is practically inapplicable on nonisolated
solutions of unregularized nonlinear systems.  
With a simple modification, a version of Newton's iteration regains its 
local quadratic convergence to nonisolated zeros of smooth mappings assuming 
the solutions are semiregular as properly defined regardless of whether the 
system is square, underdetermined or overdetermined.  
Furthermore, the iteration serves as a de facto regularization mechanism for 
computing singular zeros from empirical data. 
Even if the given system is perturbed so that the nonisolated solution 
disappears, the iteration still locally converges to a stationary point that 
approximates a solution of the underlying system with an error bound in the 
same order of the data accuracy.  
Geometrically, the iteration approximately converges to the nearest point 
on the solution manifold.  
This extension simplifies nonlinear system modeling by eliminating the zero 
isolation process and enables a wide range of applications in algebraic 
computation.
\end{abstract}

\section{Introduction}\label{s:i}

Newton's iteration as we know it loses its quadratic rate of convergence, if
it applies and converges at all, to {\em nonisolated} solutions of nonlinear 
systems of equations with large errors in numerical computation.
A subtle tweak of its formulation restores the fast rate of convergence 
and the optimal accuracy as we shall elaborate in this paper.

Perhaps there is no need to explain or even mention the importance of 
Newton's iteration in scientific computing, applied mathematics and numerical 
analysis.
In its most common and well-known formulation, Newton's iteration
(see, e.g. \cite{OrtRhe} and most textbooks in numerical analysis)
\begin{equation}\label{ni0}
\bdx_{k+1} ~~=~~ \bdx_k - J(\bdx_k)^{-1}\,\bdf(\bdx_k) \;\;\;\;\mbox{for}\;\;\;\;
k\,=\,0,1,\ldots
\end{equation}
is the standard method for solving systems of nonlinear equations 
in the form of $\bdf(\bdx)\,=\,\bdo$ where $J(\bdx)$ is the Jacobian 
of the mapping $\bdf$ at $\bdx$.
It is well documented that Newton's iteration 
quadratically converges to any isolated solution under natural 
conditions: The mapping is smooth and the initial iterate is near a solution 
at which the Jacobian is invertible.

Solving systems of nonlinear equations is a standard topic in textbooks of
numerical analysis but the discussions have always been limited to isolated
solutions of square systems. 
Models with nonisolated solutions frequently arise in 
scientific computing as we shall elaborate in \S\ref{s:m} with case studies.
However, the version (\ref{ni0}) is formulated under the assumption that 
the Jacobian is invertible at the solution and not intended for nonisolated 
solutions at which the inverse of the Jacobian is either undefined or 
nonexistent. 
Even if it converges, Newton's iteration (\ref{ni0}) is known to 
approach nonisolated solutions slowly at linear rate with an
attainable accuracy being limited and often dismal.
To circumvent those difficulties, scientific computing practitioners go to
great lengths to isolate solutions with auxiliary equations and variables. 
There have been attempts and results in the literature extending Newton's 
iteration on nonisolated solutions directly.
Those works appears to be under disseminated, scarcely applied and needing
further development.
Filling the analytical and algorithmic gap is long overdue in direct 
computation of nonisolated solutions of nonlinear systems. 

In this paper, we formulate the notion of {\em semiregular} 
solutions, establish an extension of Newton's iteration for such solutions 
and prove its local quadratic convergence on exact equations along with 
local linear convergence on perturbed equations with empirical data. 
Furthermore, we provide a geometric interpretation of the convergence 
tendency, elaborate the modeling and applications involving nonisolated 
solutions and demonstrate our software implementation with computing examples. 

An isolated zero is regular if the Jacobian is invertible, namely the nullity
of the Jacobian and the dimension of the zero are both zero. 
Nonisolated zeros of a smooth mapping can form a smooth submanifold of
a positive dimension such as curves and surfaces.
We generalize the regularity of isolated zeros to the semiregularity of 
nonisolated cases as the dimension being identical to the nullity of the 
Jacobian at the zero.
Semiregular zeros of a positive dimension form branches with locally 
invariant dimensions (c.f. Lemma~\ref{l:cr}) and,
near a semiregular zero, we prove a crucial property that every stationary 
point is a semiregular zero (c.f. Lemma~\ref{l:rz}). 

We extended Newton's iteration to the form of
\begin{equation}\label{ni1}
\bdx_{k+1} ~~=~~ \bdx_k - J_\rkr(\bdx_k)^\dagger\,\bdf(\bdx_k) 
\;\;\;\;\mbox{for}\;\;\;\; k\,=\,0,1,\ldots
\end{equation}
for a smooth mapping $\bdf$ from an open domain in $\R^m$ or $\C^m$ 
to $\R^n$ or $\C^n$ where $J_\rkr(\bdx_k)^\dagger$ is the 
Moore-Penrose inverse of $J_\rkr(\bdx_k)$ that is the rank-$r$ projection 
of the Jacobian $J(\bdx_k)$. 
We establish its local quadratic convergence (Theorem~\ref{t:mt}) under
minimal natural assumptions: The mapping $\bdf$ is smooth and the initial 
iterate is near a semiregular zero at which the Jacobian is of rank $r$. 

Nonisolated solutions can be highly sensitive to data perturbations.
When the system of equations is perturbed, represented with empirical data
or solved using floating point arithmetic, the nonisolated solution can be
significantly altered or even disappears altogether. 
We prove that the proposed Newton's iteration still converges to a stationary
point that approximates an exact solution of the underlying system with an
accuracy in the same order of the data error (c.f. Theorem~\ref{t:mt2}). 
In other words, the proposed extension of Newton's iteration also serves as
a regularization mechanism for such an ill-posed zero-finding problem.
A condition number can also be derived from Theorem~\ref{t:mt2} for 
nonisolated zeros with respect to data perturbations. 

As a geometric interpretation, we shall illustrate the behavior of the
Newton's iteration (\ref{ni1}) for converging to roughly the point on the
solution manifold nearest to the initial iterate.
We also elaborate case studies on mathematical modeling with nonisolated
solutions, demonstrate our software implementation in solving for such
solutions with step-by-step calling sequences and computing results.

Extending Newton's iteration beyond (\ref{ni0}) by replacing the inverse
with a certain kind of generalized inverses traces back to Gauss for solving
least squares solutions of overdetermined systems by the 
Gauss-Newton iteration with an assumption that the Jacobian is injective.
For systems with rank-deficient Jacobians, Ben-Israel \cite{Ben66} is the 
first to propose using Moore-Penrose inverses to generalize Newton's iteration 
as
\begin{equation}\label{bi}
\bdx_{k+1} ~~=~~ \bdx_k - J(\bdx_k)^\dagger\,\bdf(\bdx_k) \;\;\;\;\mbox{for}
\;\;\;\; k\,=\,0,1,\ldots
\end{equation}
``but the conditions for the [convergence] theorem are somewhat restrictive
and unnatural'' \cite{Boggs76}.
Chu is the first to prove the local convergence of Newton's iteration
(\ref{bi}) with essentially minimal assumptions for 
underdetermined systems with surjective Jacobians \cite{Chu83}. 
Applying the alpha theory, Dedieu and Kim \cite{DedKim} prove Newton's 
iteration (\ref{bi}) locally quadratically converges to a nonisolated 
solution under the assumption that the Jacobian has a constant
deficient rank in a neighborhood of the initial iterate.
 
~Nashed and Chen \cite{NasChe93} propose a Newton-like iteration 
by replacing the inverse with a certain outer inverse of
the Jacobian and prove a quadratic convergence to a stationary point
under a set of conditions,
but it is unknown whether the stationary point is a
zero of ~$\bdf$.
~Chen, Nashed and Qi \cite{CheNasQi} along with Levin and Ben-Israel 
\cite{LevBen} follow up with similar convergence 
results toward stationary points using outer inverses.

In comparison to those pioneer works, our extension (\ref{ni1}) is suitable
for any rank of the Jacobian at the zero. 
Furthermore, our convergence theorems require minimal assumptions and 
the iteration quadratically converges to a stationary point that is 
guaranteed to be a zero if the initial iterate is near a semiregular zero.
The iteration also permits data perturbations and floating point
arithmetic while still converges to an approximate zero at linear rate.

We loosely refer to our extension (\ref{ni1}) as the
{\em rank-$r$ Newton's iteration} or simply {\em Newton's iteration} following 
a long-standing practice.
The terminology is actually debatable as it is fair to ask
``Is Newton's method really Newton's method?''\,\cite{Deu12}, 
and the term may even be considered ``an enduring myth''
\cite{Kol92}.
Heavily influenced by Fran\c{c}ois Vi\`ete, Isaac Newton's original method 
is not even an iteration and can be considered a special case of 
Joseph Raphson's formulation restricted to univariate polynomial equations. 
The version (\ref{ni0}) of Newton's iteration can deservedly be credit to 
Thomas Simpson as well.
As a myth or not, however, the term ``Newton's iteration'' has been 
used for all extensions of Newton's original method and will likely to used 
as a convention in the future.  

\section{Preliminaries}\label{s:p}

Column vectors are denoted by boldface lower case letters such as $\bdb$,
$\bdx$, $\bdy$ etc. with $\bdo$ being a zero vector whose dimension can
be derived from the context.
The vector spaces of $n$-dimensional real and complex column vectors are
denoted by $\R^n$ and $\C^n$ respectively.
The vector space of $m\times n$ complex matrices including real matrices
is denoted by $\C^{m\times n}$.
Matrices are denoted by upper case letters such as $A$, $B$, $X$, etc.
with $O$ and $I$ denoting a zero matrix and an identity matrix 
respectively.
The range, kernel, rank and Hermitian transpose of a matrix $A$ are
denoted by $\cR(A)$, $\cK(A)$, $\rank{A}$ and $A^\h$ respectively.
For any matrix $A$, its {\em Moore-Penrose inverse}
\cite[\S 5.5.2, p. 290]{golub-vanloan4}
$A^\dagger$ is the unique matrix satisfying
\begin{equation}\label{mpc}
A\,A^\dagger\,A ~=~ A, ~~
A^\dagger\,A\,A^\dagger ~=~ A^\dagger, ~~
(A\,A^\dagger)^\h ~=~ A\,A^\dagger, ~~ 
(A^\dagger\,A)^\h ~=~ A^\dagger\,A.
\end{equation}
The \,$j$-th largest singular value of $A$ is denoted by $\sg_j(A)$.
Let $U\,\Sigma\,V^\h$ be the singular value decomposition of $A$
where $U\,=\,[\bdu_1,\,\cdots,\,\bdu_m]$ and 
$V\,=\,[\bdv_1,\,\cdots,\,\bdv_n]$ are unitary matrices formed by 
the left singular vectors and the right singular vectors respectively. 
The {\em rank-$r$ projection} $A_\rkr$ of $A$, also known as rank-$r$
truncated singular value decomposition (TSVD) and rank-$r$ approximation of 
$A$, is defined as
\[ A_\rkr ~~:=~~ \sg_1(A)\,\bdu_1\,\bdv_1^\h + \cdots +
\sg_r(A)\,\bdu_r\,\bdv_r^\h 
\]
Using singular values and singular vectors, the identity
~\cite[\S 5.5.2]{golub-vanloan4}
\[  A^\dagger ~~\equiv~~ 
\sum_{\sg_j(A)>0} \,\frac{1}{\sg_j(A)}\,\bdv_j\,\bdu_j^\h
\]
holds and it is straightforward to verify
\begin{align}\label{ara}
A\,A_\rkr^\dagger &~~=~~ A_\rkr\,A_\rkr^\dagger ~~=~~ 
[\bdu_1,\cdots,\bdu_r]\, [\bdu_1,\cdots,\bdu_r]^\h \\
A_\rkr^\dagger\,A &~~=~~ A_\rkr^\dagger\,A_\rkr ~~=~~ 
[\bdv_1,\cdots,\bdv_r]\, [\bdv_1,\cdots,\bdv_r]^\h \label{ara2}
\end{align}
that are orthogonal projections onto the subspaces spanned by  
$\{\bdu_1,\ldots,\bdu_r\}$ and
$\{\bdv_1,\ldots,\bdv_r\}$ respectively.
For any matrix $A$, denote $A_\rkr^\dagger$ as $(A_\rkr)^\dagger$.

We say $\bdf$ is a {\em smooth mapping} \,if 
$\bdf\,:\,\Omega\subset\R^n\, \rightarrow\,\R^m$ has continuous derivatives 
of second order, or $\bdf\,:\,\Omega\subset\C^n\, \rightarrow\,\C^m$ 
is holomorphic, where the domain $\Omega$ is an open subset of $\C^n$ or
$\R^n$.
We may designate a variable name, say $\bdx$, for $\bdf$ 
and denote $\bdf$ as $\bdx\,\mapsto\,\bdf(\bdx)$.
In that case the Jacobian of $\bdf$ at any particular 
$\bdx_0\,\in\,\Omega$ is denoted by $\bdf_\bdx(\bdx_0)$ or $J(\bdx_0)$
while $J_\rkr(\bdx_0)$ or equivalently $\bdf_\bdx(\bdx_0)_\rkr$ is 
its rank-$r$ projection.
For a smooth mapping $(\bdx,\bdy)\,\mapsto\,\bdf(\bdx,\bdy)$ at 
$(\bdx_0,\bdy_0)$, the notation $\bdf_{\bdx\bdy}(\bdx_0,\bdy_0)$ 
represents its Jacobian (with respect to {\em both} $\bdx$ and $\bdy$) while
$\bdf_\bdx(\bdx_0,\bdy_0)$ and $\bdf_\bdy(\bdx_0,\bdy_0)$ denote the 
(partial) Jacobians with respect to $\bdx$ and $\bdy$ respectively at
$(\bdx_0,\bdy_0)$.

\begin{lemma}\label{l:p}
Let ~$J(\bdx)$ ~be the Jacobian of a smooth mapping ~$\bdf$ ~at any
$\bdx$ ~in its open domain ~$\Omega$ ~in ~$\C^n$ ~or ~$\R^n$.
Assume ~$\rank{J(\bdx_*)}=r$ ~at ~$\bdx_*\in\Omega$.
Then there is an open bounded convex subset ~$\Omega_*\ni\bdx_*$ ~of 
~$\Omega$ ~and constants ~$\zeta,\,\mu,\,\eta,\,\xi>0$ ~such that,
for every ~$\bdx,\,\bdy\in\Omega_*$, ~the inequality
$\rank{J(\bdx)} \ge r$ ~holds along with
\begin{align}
  \big\|J_\rkr(\bdx)\,J_\rkr(\bdx)^\dagger-J_\rkr(\bdy)\,J_\rkr(\bdy)^\dagger
\big\|_2
&~~\le~~ \zeta\,\|\bdx-\bdy\|_2 \label{jrxy} \\
  \big\|J_\rkr(\bdx) -J_\rkr(\bdy)\big\|_2 
&~~\le~~ \mu\,\|\bdx-\bdy\|_2 \label{jrxy1} \\
  \big\|J_\rkr(\bdx)^\dagger -J_\rkr(\bdy)^\dagger\big\|_2 
&~~\le~~ \eta\,\|\bdx-\bdy\|_2 \label{jrxy2} \\
  \big\|J_\rkr(\bdx)^\dagger \big\|_2
&~~\le~~ \xi \label{jrx3}
\end{align}
\end{lemma}

{\em Proof.} 
From ~$\rank{J(\bdx_*)}=r$, ~we have 
$\sg_r(J(\bdx_*))>\sg_{r+1}(J(\bdx_*))=0$. 
Weyl's Theorem \cite[Corollary 4.31, p. 69]{stew1} ensures the singular
values to be continuous with respect to the matrix entries. 
By the continuity of ~$J(\bdx)$ ~with respect to ~$\bdx$ ~in
$\Omega$, ~there is an open bounded convex neighborhood ~$\Omega_*$ ~of 
$\bdx_*$ ~with ~$\overline{\Omega}_*\subset\Omega$ ~such
that ~$\sg_r(J(\bdx))>2\sg_{r+1}(J(\bdx))$ ~for every 
~$\bdx\in\Omega_*$.
We can further assume ~$\Omega_*$ ~to be sufficiently small so that
\begin{equation}\label{wc}  \big\|J(\bdx)-J(\bdy)\big\|_2 ~~<~~
\mbox{$\frac{1}{2}$} \big(\sg_r(J(\bdx))-\sg_{r+1}(J(\bdx))\big)
\end{equation}
for all ~$\bdx,\,\bdy\,\in\,\Omega_*$ ~and
\begin{equation}\label{mj}
\max_{\bdx\,\in\,\overline{\Omega}_*} \big\|J_\rkr(\bdx)^\dagger\big\|_2
~~=~~ 
\max_{\bdx\,\in\,\overline{\Omega}_*} 
\frac{1}{\sg_r(J(\bdx))} ~~\le~~ 2\,\big\|J(\bdx_*)^\dagger\big\|_2
\end{equation}
so (\ref{jrx3}) is true.
The left-hand side of (\ref{jrxy}) is the distance between the subspaces
spanned by the first \,$r$ \,left singular vectors of ~$J(\bdx)$ ~and 
~$J(\bdy)$ ~respectively and, by Wedin's Theorem \cite{wedin} (also see
\cite[Theorem 4.4]{StewSun}) and (\ref{wc}),  
there is a constant ~$\zeta\,>\,0$ ~such that
\begin{align*}
\big\|J_\rkr(\bdx)\,J_\rkr(\bdx)^\dagger &-J_\rkr(\bdy)\,J_\rkr(\bdy)^\dagger
\big\|_2\\
&~~\le~~ 
4\,\big\|J_\rkr(\bdx)^\dagger\big\|_2\,\big\|J(\bdx)-J(\bdy)\big\|_2 \\
&~~\le~~ \zeta\,\|\bdx-\bdy\|_2
\end{align*}
for all ~$\bdx,\,\bdy\,\in\,\Omega_*$. 
As a result, the inequality (\ref{jrxy1}) follows from (\ref{jrxy}) and,
by (\ref{mpc}) and (\ref{ara2}),
\begin{align*}
\lefteqn{\big\|J_\rkr(\bdx) -J_\rkr(\bdy)\big\|_2}\\
&~~= ~~ \big\|J_\rkr(\bdx)\,J_\rkr(\bdx)^\dagger\,J(\bdx)
-J_\rkr(\bdy)\,J_\rkr(\bdy)^\dagger \,J(\bdy)\big\|_2 
\\
&~~\le~~ \big\|J_\rkr(\bdx)\,J_\rkr(\bdx)^\dagger\big\|_2\,
\|J(\bdx)-J(\bdy)\|_2 \\
&~~~~~~~~~ +
\big\|J_\rkr(\bdx)\,J_\rkr(\bdx)^\dagger-J_\rkr(\bdy)\,J_\rkr(\bdy)^\dagger
\big\|_2\,\|J(\bdy)\|_2. 
\end{align*}
since ~$\big\|J_\rkr(\bdx)\,J_\rkr(\bdx)^\dagger\big\|_2\,=\,1$ ~and
$\|J(\bdy)\|_2$ ~is bounded on the compact set ~$\overline{\Omega}_*$.
By \cite[Theorem 3.3]{stew77}, there is a constant ~$\al\,>\,0$ ~such that
\begin{align*}
\lefteqn{\big\|J_\rkr(\bdx)^\dagger -J_\rkr(\bdy)^\dagger\big\|_2} \\
&~~\le~~ \al\,\big\|J_\rkr(\bdx)^\dagger\big\|_2\,\big\|J_\rkr(\bdy)^\dagger
\big\|_2 \,\big\|J_\rkr(\bdx)-J_\rkr(\bdy)\big\|_2
\end{align*}
leading to (\ref{jrxy2}).
~\qed 

\begin{lemma}\label{l:cp}
For $\F\,=\,\C$ or $\R$, let $\bdz\,\mapsto\,\phi(\bdz)$ be 
a continuous injective mapping from an open set $\Omega$ in $\F^n$ to 
$\F^m$.
At any $\bdz_0\,\in\,\Omega$, 
there is an open neighborhood $\Dl$ of 
$\phi(\bdz_0)$ in $\F^m$ such that, for every 
$\bdb\,\in\,\Dl$, there exists a $\bdz_\bdb$ in $\Omega$ and an
open neighborhood $\Sigma_0$ of $\bdz_\bdb$ with
\begin{equation}\label{bp}
  \|\bdb-\phi(\bdz_\bdb)\|_2 ~~=~~ \min_{\bdz\in\Sigma_0}
\|\bdb-\phi(\bdz)\|_2.
\end{equation}
Further assume $\phi$ is differentiable in $\Omega$. 
Then 
\begin{equation}\label{cpc}
\phi_\bdz(\bdz_\bdb)^\dagger\,(\bdb-\phi(\bdz_\bdb))\,=\,\bdo.
\end{equation}
\end{lemma}

{\em Proof.}
Let $\Sigma_0$ be an open bounded neighborhood of 
$\bdz_0$ such than $\overline{\Sigma}_0\,\subset\,\Omega$.
Since $\phi$ ~is one-to-one and continuous, we have
\[ \dl ~~=~~ \min_{\bdz\in\overline{\Sigma}_0\setminus\Sigma_0}\,
\|\phi(\bdz)-\phi(\bdz_0)\|_2 ~~>~~ 0.
\]
Let $\Dl\,=\,\big\{\bdy\,\in\,\F^m ~\big|~ \|\bdy-\phi(\bdz_0)\|_2\,<\,
\frac{1}{2}\,\dl\big\}$.
Then, for every $\bdb\,\in\,\Dl$, there exists a $\bdz_\bdb\,\in\,
\overline{\Sigma}_0$ such that
$\|\phi(\bdz_\bdb)-\bdb\|_2 \,=\, \min_{\bdz\in\overline{\Sigma}_0}\,
\|\phi(\bdz)-\bdb\|_2$.
For every $\bdz\,\in\,\overline{\Sigma}_0\setminus\Sigma_0$, however,
\begin{align*}   \|\phi(\bdz)-\bdb\|_2 &~~\ge~~ \|\phi(\bdz)-\phi(\bdz_0)\|_2-
\|\phi(\bdz_0)-\bdb\|_2  \\
&~~>~~ \mbox{$\frac{1}{2}$}\,\dl
~~>~~ \|\phi(\bdz_0)-\bdb\|_2 ~~\ge~~  \|\phi(\bdz_\bdb)-\bdb\|_2 
\end{align*}
implying $\bdz_\bdb\,\in\,\Sigma_0$ and (\ref{bp}).
Since a local minimum of 
\[ \|\phi(\bdz)-\bdb\|_2^2 ~~=~~  (\phi(\bdz)-\bdb)^\h\, (\phi(\bdz)-\bdb)
\]
occurs at the interior point $\bdz_\bdb\,\in\,\Sigma_0$, it is 
straightforward to verify the equation
$\phi_\bdz(\bdz_\bdb)^\h\,(\phi(\bdz_\bdb)-\bdb)\,= \,\bdo$ and 
thus (\ref{cpc}) 
from $\cR(\phi_\bdz(\bdz_\bdb)^\h)\,= \,\cR(\phi_\bdz(\bdz_\bdb)^\dagger)$.
~\qed

\section{Semiregular zeros of smooth mappings}\label{s:r}

A point $\bdx_*$ is a {\em zero} of a mapping $\bdf$ if 
$\bdx\,=\,\bdx_*$ is a {\em solution} of the equation
$\bdf(\bdx)\,=\,\bdo$.
A common notation $\bdf^{-1}(\bdo)$ stands for the set of all zeros of 
$\bdf$.
A zero $\bdx_*$ of $\bdf$ is {\em isolated} if there is an open 
neighborhood $\Lambda$ of $\bdx_*$ in the domain of $\bdf$ such that
$\bdf^{-1}(\bdo)\cap\Lambda\,=\,\{\bdx_*\}$ or $\bdx_*$
is {\em nonisolated} otherwise.
A nonisolated zero $\bdx_*$ of a smooth mapping $\bdf$ may belong to
a curve, a surface or a higher dimensional subset of $\bdf^{-1}(\bdo)$. 
We adopt a simple definition of the dimension of a zero as follows.
For more in-depth elaboration on the dimension of zero sets, see 
\cite[p. 17]{BSHW13}.

\begin{definition}[Dimension of a Zero]\label{d:dz}
For $\F\,=\C$ or $\R$, let $\bdx_*$ be a zero of a
smooth mapping $\bdf\,:\,\Omega\subset\F^m\,\rightarrow\,\F^n$.
If there is an open neighborhood $\Dl\subset\Omega$ of $\bdx_*$ in 
$\F^m$ such that $\Dl\cap\bdf^{-1}(\bdo)\,=\,\phi(\Lambda)$ where 
$\bdz\,\mapsto\,\phi(\bdz)$ is a differentiable injective mapping defined 
in a connected open set $\Lambda$ in $\F^k$ for a certain $k\,>\,0$ with 
$\phi(\bdz_*)\,=\,\bdx_*$ and $\rank{\phi_\bdz(\bdz_*)}\,=\,k$,
then the {\em dimension} of $\bdx_*$ as a zero of $\bdf$ is defined
as
\[
  \dm_\bdf(\bdx_*) ~~:=~~ \dm(\cR(\phi_\bdz(\bdz_*))) ~~\equiv~~
\rank{\phi_\bdz(\bdz_*)} ~~=~~ k.
\]
As a special case, an isolated zero is of dimension zero.
\end{definition}

A so-defined $k$-dimensional zero $\bdx_*$ is on a smooth submanifold of
dimension $k$ in $\F^m$.
If the dimension $\dm_\bdf(\bdx_*)$ is well-defined, then 
$\phi_\bdz(\bdz_*)$ in Definition~\ref{d:dz} is of full column rank and 
there is an open neighborhood $\Lambda_*\subset\Lambda$ of $\bdz_*$ 
such that 
$\rank{\phi_\bdz(\hat\bdz)}\,\equiv\,k$ for all $\hat\bdz\in\Lambda_*$.
Namely {\em the dimension of a zero is locally invariant}.
We shall also say every 
$\bdx\,\in\,\phi(\Lambda_*)$ is {\em in the same branch} of zeros as 
$\bdx_*$ and, if a zero $\tilde\bdx$ is in the same branch of $\bdx_*$, 
every zero $\hat\bdx$ in the same branch of $\tilde\bdx$ is in the same 
branch of $\bdx_*$.

An isolated zero $\bdx_*$ of $\bdf$ is regular if its dimension $0$ is 
identical to the nullity of the Jacobian $\bdf_\bdx(\bdx_*)$.
This characteristic of regularity can naturally be generalized to zeros of 
higher dimensions for being semiregular as defined below.
There are tremendous advantages in computing semiregular zeros as we shall 
elaborate throughout this paper. 

\begin{definition}[Semiregular Zero]\label{d:rz}
A zero $\bdx_*$ of a smooth mapping $\bdx\,\mapsto\,\bdf(\bdx)$
is {\em semiregular} if $\dm_\bdf(\bdx_*)$ is well-defined and
identical to $\nullity{\bdf_\bdx(\bdx_*)}$. 
Namely
\begin{equation}\label{drn}
 \dm_\bdf(\bdx_*)+\rank{\bdf_\bdx(\bdx_*)} ~~=~~ 
\mbox{the dimension of the domain of ~$\bdf$}.
\end{equation}
A zero is {\em ultrasingular} if it is not semiregular.
\end{definition}

A system of equations $\bdf(\bdx)\,=\,\bdo$ is said to be 
{\em underdetermined} \,if $\bdf$ is a mapping from $\Omega\subset\F^m$ to
$\F^n$ with $m\,>\,n$ where $\F\,=\,\C$ or $\R$.
A solution of an underdetermined system is always semiregular if the Jacobian 
is surjective or, equivalently, of full row rank.
For instance, let $(\bdu_*,\bdv_*)$ be a zero of a smooth mapping 
$(\bdu,\bdv)\,\mapsto\,\bdf(\bdu,\bdv)$ from $\R^k\times\R^m$
to $\R^m$ and the partial Jacobian $\bdf_\bdv(\bdu_*,\bdv_*)$ is 
invertible. 
By the Implicit Mapping Theorem, there is a differentiable mapping 
$\bdu\,\mapsto\,\bdg(\bdu)$ from a neighborhood $\Lambda$ of 
$\bdu_*$ in $\R^k$ to $\R^m$ with $\bdg(\bdu_*)\,=\,\bdv_*$ and 
there is an open neighborhood $\Dl$ of $(\bdu_*,\bdv_*)$ in 
$\R^k\times\R^m$ such that $\Dl\cap\bdf^{-1}(\bdo)\,=\,\phi(\Lambda)$ where 
$\phi(\bdu)\,=\,(\bdu,\bdg(\bdu))$ for $\bdu\in\Lambda$.
Furthermore, the Jacobian $\phi_\bdu(\bdu_*)$ is obviously of full 
column rank $k$.
As a result, the dimension of the zero $(\bdu_*,\bdv_*)$ is $k$ that
is identical to the nullity of Jacobian of $\bdf$ at $(\bdu_*,\bdv_*)$,
implying $(\bdu_*,\bdv_*)$ is semiregular.

Ultrasingular zeros can be isolated multiple zeros \cite{DLZ}, isolated 
ultrasingularity embedded in nonisolated zero set (c.f. Example~\ref{e:is}) or 
can form an entire branch of zeros (c.f. Example~\ref{e:sb}).
Like the dimension, semiregularity is also invariant on a branch of
nonisolated zeros as asserted in the following lemma.

\begin{lemma}[Local Invariance of Semiregularity]\label{l:cr}
Let $\bdx_*$ be a semiregular zero of a smooth mapping $\bdf$.
Then there is an open neighborhood $\Dl_*$ of $\bdx_*$ such 
that every $\hat\bdx\,\in\,\Dl_*\cap\bdf^{-1}(\bdo)$ is a semiregular zero of 
$\bdf$ in the same branch of $\bdx_*$.
\end{lemma}

{\em Proof}. 
The assertion is obviously true for isolated regular zeros.
Let $\F$ be either $\C$ \,or\, $\R$ and the domain of $\bdf$ is
an open subset in $\F^n$.
Assume $\bdx_*$ is a semiregular $k$-dimensional zero of $\bdf$.
There is an open neighborhoods $\Dl$ of $\bdx_*$ in $\F^n$
and there is an open connected set $\Lambda_1$ of a certain $\bdz_*$ in 
$\F^k$ along with a differentiable injective mapping 
$\phi\,:\,\Lambda_1\,\rightarrow\,\F^n$ such that
$\phi(\bdz_*)\,=\,\bdx_*$, $\Dl\cap\bdf^{-1}(\bdo)\,=\,\phi(\Lambda_1)$,
$\rank{\bdf_\bdx(\bdx_*)}\,=\,n-k$ and $\rank{\phi_\bdz(\bdz_*)}\,=\,k$.
By the continuity of singular values we can assume 
$\rank{\phi_\bdz(\bdz)}\,\equiv\,k$ for all $\bdz\in\Lambda_1$.
From $\bdf(\phi(\bdz))\,\equiv\,\bdo$ and 
$\bdf_\bdx(\phi(\bdz))\,\phi_\bdz(\bdz)\,\equiv\,O$ for all $\bdz\in
\Lambda_1$, we have $\nullity{\bdf_\bdx(\phi(\bdz))}\,\ge\,k$ for all
$\bdz\in\Lambda_1$.
By the continuity of singular values again, there is an open neighborhood
$\Dl_*\subset\Dl$ of $\bdx_*$ such that 
$\rank{\bdf_\bdx(\bdx)}\,\ge\,n-k$ for all $\bdx\in\Dl_*$.
Consequently, every $\hat\bdx\in\Lambda\cap\bdf^{-1}(\bdo)$ is a semiregular 
zero where $\Lambda\,=\,\phi^{-1}(\Dl_*)$.
~\qed

We shall propose a new version of Newton's iteration 
that, under proper conditions, converges
to a {\em stationary point} $\hat\bdx$ at which 
$J_\rkr(\hat\bdx)^\dagger\,\bdf(\hat\bdx)\,=\,\bdo$.
The following stationary point property of semiregular zeros ensures that, 
in a neighborhood of a semiregular zero, all stationary points are 
semiregular zeros in the same branch.

\begin{lemma}[Stationary Point Property]\label{l:rz}
Let $\bdx\,\mapsto\,\bdf(\bdx)$ be a smooth mapping with a semiregular
zero $\bdx_*$ and $r\,=\,\rank{\bdf_\bdx(\bdx_*)}$.
Then there is an open neighborhood $\Omega_*$ of $\bdx_*$ such that,
for  any $\hat\bdx\,\in\,\Omega_*$, the equality
$\bdf_\bdx(\hat\bdx)_\rkr^\dagger\,\bdf(\hat\bdx)\,=\,\bdo$ 
holds if and only if $\hat\bdx$ is a semiregular zero of $\bdf$ in the
same branch of $\bdx_*$. 
\end{lemma}

{\em Proof.}
We claim there is a neighborhood $\Omega_1$ of $\bdx_*$ such
that $\bdf(\hat\bdx)\,=\,\bdo$ for every $\hat\bdx\in\Omega_1$ with
$\bdf_\bdx(\hat\bdx)_\rkr^\dagger\,\bdf(\hat\bdx)\,=\,\bdo$.
Assume this assertion is false. 
Namely there is a sequence 
$\{\bdx_j\}_{j=1}^\infty$ converging to $\bdx_*$ such that
$\bdf_\bdx(\bdx_j)_\rkr^\dagger\,\bdf(\bdx_j)\,=\,\bdo$ but 
$\bdf(\bdx_j)\,\ne\,\bdo$ for all $j\,=\,1,2,\ldots$.
Let $\bdz\,\mapsto\,\phi(\bdz)$ be the parameterization of 
the solution branch containing $\bdx_*$ as in Definition~\ref{d:dz}
with $\phi(\bdz_*)\,=\,\bdx_*$.
%
%
From Lemma~\ref{l:cp} with any sufficiently large $j$, 
there is a $\check\bdx_j\,\in\,\Omega_*\cap \bdf^{-1}(\bdo)\,=\,\phi(\Dl)$ 
such that 
\begin{equation}\label{jj2}   \|\bdx_j-\check\bdx_j\|_2 ~~=~~ 
\min_{\bdz\in\Dl}\|\bdx_j-\phi(\bdz)\|_2 ~~=~~ \|\bdx_j-\phi(\bdz_j)\|_2
\end{equation}
at a certain $\bdz_j$ with $\phi(\bdz_j)\,=\,\check\bdx_j$, implying 
\begin{equation}\label{pxp}
\phi_\bdz(\bdz_j)\,\phi_\bdz(\bdz_j)^\dagger\,
\frac{\bdx_j-\phi(\bdz_j)}{\|\bdx_j-\phi(\bdz_j)\|_2} ~~=~~
\frac{\phi_\bdz(\bdz_j)}{\|\bdx_j-\phi(\bdz_j)\|_2} 
\,\Big(\phi_\bdz(\bdz_j)^\dagger\,\big(\bdx_j-\phi(\bdz_j)\big)\Big)
~~=~~ \bdo.
\end{equation}
We claim $\check\bdx_j\,\rightarrow\,\bdx_*$ as well when $j\,\rightarrow
\,\infty$.
Assume otherwise. 
Namely there is an $\eps\,>\,0$ such that, for any $N\,>\,0$, there 
is a $j\,>\,N$ with $\|\check\bdx_j-\bdx_*\|_2\,\ge\,2\,\eps$.
However, we have $\|\bdx_j-\bdx_*\|_2\,<\,\eps$ for all
$j$ larger than a certain $N$, implying 
\[ \|\check\bdx_j-\bdx_j\| ~~\ge~~ \|\check\bdx_j-\bdx_*\| -
\|\bdx_*-\bdx_j\| ~~>~~ \eps ~~>~~ \|\bdx_j-\bdx_*\|_2
\]
that is a contradiction to (\ref{jj2}).

Since $\bdf(\bdx_j)\,\ne\,\bdo$, we have $\bdx_j\,\ne\,\check\bdx_j$ and
we can assume $(\bdx_j-\check\bdx_j)/\|\bdx_j-\check\bdx_j\|_2$
converges to a unit vector $\bdv$ for $j\,\rightarrow\,\infty$ due to
compactness.
Then
\begin{align} 
0&~~=~~  \lim_{j\rightarrow\infty}\,
\frac{\bdf_\bdx(\bdx_j)_\rkr^\dagger\,\big(\bdf(\check\bdx_j)
-\bdf(\bdx_j)\big)}{ \|\bdx_j-\check\bdx_j\|_2} \nonumber \\
& ~~=~~ \lim_{j\rightarrow\infty}\,\frac{\bdf_\bdx(\bdx_j)_\rkr^\dagger\,
\bdf_\bdx(\bdx_j) \,(\check\bdx_j-\bdx_j)}{\|\check\bdx_j-\bdx_j\|_2} 
~~=~~ \bdf_\bdx(\bdx_*)^\dagger\,\bdf_\bdx(\bdx_*)\,\bdv \label{fxfxv}
\end{align}
by (\ref{ara2}) and (\ref{jrxy}), implying 
$\bdv\,\in\,\cK(\bdf_\bdx(\bdx_*))$.
As a result, 
\[  \spn\{\bdv\}\oplus\cR(\phi_\bdz(\bdz_*)) ~~\subset~~
\cK(\bdf_\bdx(\bdx_*))
\]
since $\bdf_\bdx(\bdx_*)\phi_\bdz(\bdz_*)\,=\,O$ due to 
$\bdf(\phi(\bdz))\,\equiv\, \bdo$ in a neighborhood of $\bdz_*$.
From the limit of (\ref{pxp}) for $j\,\rightarrow\,\infty$, we have 
$\bdv\,\in\,\cR(\phi_\bdz(\bdz_*))^\perp$ and thus
\[  \nullity{\bdf_\bdx(\bdx_*)} ~~\ge~~ \rank{\phi_\bdz(\bdz_*)}+1
\]
which is a contradiction to the semiregularity of $\bdx_*$. 

For the special case where $\bdx_*$ is an isolated semiregular zero of 
$\bdf$ with dimension 0, the above proof applies
with $\check\bdx_j\,=\,\bdx_j$ so that (\ref{fxfxv}) holds, implying a
contradiction to $\nullity{\bdf_\bdx(\bdx_*)}\,=\,0$.

By Lemma~\ref{l:cr}, there is a neighborhood $\Omega_2$ of $\bdx_*$ such
that every $\bdx\in\Omega_2\cap\bdf^{-1}(\bdo)$ is a semiregular zero of 
$\bdf$ in the same branch of $\bdx_*$. 
Thus the lemma holds for $\Omega_*\,=\,\Omega_1\cap\Omega_2$.
~\qed

\begin{remark}[A note on terminology]\em
Keller \cite{Kel81} argues that the terminology ``is somewhat unfortunate''
on {\em isolated} and {\em nonisolated} zeros as he defines the later as
zeros at which the Jacobian is not injective. 
The isolated zeros here are referred to as {\em geometrically isolated} in
\cite{Kel81}.
{\em Nonisolated}\, zeros defined by Keller include multiple 0-dimensional
zeros and zeros on positive dimensional branches.
The term {\em singular zero} is broadly accepted to include multiple isolated
zero and nonisolated zero (c.f. \cite[\S 1.2.2, p. 20]{BSHW13}).
We propose the term {\em ultrasingular zero} to distinguish those singular 
zeros from semiregular ones as defined in Definition~\ref{d:rz}.
\end{remark}

\section{Convergence theorem on exact equations}

Consider the system of equations in the form of \;$\bdf(\bdx)\,=\,\bdo$
\;where \;$\bdf\,:\,\Omega\subset\F^m\,\rightarrow\,\F^n$ is a smooth
mapping with $\F\,=\,\C$ or $\F\,=\,\R$. 
The system can be square ($m\,=\,n$), underdetermined ($m\,>\,n$) or
overdetermined ($m\,<\,n$).
We propose the iteration
\begin{equation}\label{ni}
\bdx_{k+1} ~~=~~ \bdx_k - J_\rkr(\bdx_k)^\dagger\,\bdf(\bdx_k)
~~~~\mbox{for}~~~ k = 0, 1, \ldots
\end{equation}
for computing a zero $\bdx_*$ of $\bdf$ at which the Jacobian 
$J(\bdx_*)$ is of rank $r$ particularly when $\bdx_*$ is on a branch of 
semiregular nonisolated zeros.
We loosely refer to (\ref{ni}) as the
{\em rank-$r$ Newton's iteration} or simply
{\em Newton's iteration} since it is 
identical to the commonly-known Newton's iteration when $r\,=\,m\,=\,n$.
Assume the mapping $\bdf$ is given with exact data. 
The following theorem establishes the local quadratic convergence of the 
iteration (\ref{ni}).
We shall consider the equation with empirical data in \S\ref{s:cp}.

\begin{theorem}[Convergence Theorem]\label{t:mt}
Let $\bdf$ be a smooth mapping in an open domain 
with a rank \,$r$ \,Jacobian $J(\bdx_*)$ at a semiregular zero $\bdx_*$.
For every open neighborhood $\Omega_1$ of $\bdx_*$, 
there is a neighborhood $\Omega_0$ of $\bdx_*$ such that, 
from every initial iterate $\bdx_0\,\in\,\Omega_0$,
the rank-$r$ Newton's iteration {\em (\ref{ni})} converges 
quadratically to a semiregular zero $\hat\bdx\in\Omega_1$ of $\bdf$ 
in the same branch as $\bdx_*$.
\end{theorem}

{\em Proof.}
Let ~$\Omega_*\ni\bdx_*$ ~be an open convex subset of the domain
for ~$\bdf$ ~as specified in Lemma~\ref{l:p} so that (\ref{jrxy}) and 
(\ref{jrx3}) hold in ~$\Omega_*$. 
From Lemma~\ref{l:rz}, we can further assume 
$J_\rkr(\bdx)^\dagger\,\bdf(\bdx)=\bdo$ ~implies ~$\bdx$ ~is a semiregular
zero of ~$\bdf$ ~in the same branch of ~$\bdx_*$ ~for every ~$\bdx\in\Omega_*$.
Since ~$\bdf$ ~is smooth, there are constants ~$\mu,\,\gamma\,>\,0$ ~such
that 
\begin{align*}  
\big\|\bdf(\bdy)-\bdf(\bdx)\big\|_2 
&~~\le~~ \mu\,\|\bdy-\bdx\|_2 \\
\big\|\bdf(\bdy)-\bdf(\bdx)-J(\bdx)\,(\bdy-\bdx)\big\|_2 
&~~\le~~ \gamma\,\|\bdy-\bdx\|_2^2
\end{align*}
for all ~$\bdx,\,\bdy\,\in\,\Omega_*$.
Denote 
$S_\eps(\bdx_*) \,:=\, \big\{\bdx ~\big|~ \|\bdx-\bdx_*\|_2\,<\,\eps\big\}$
for any ~$\eps\,>\,0$.
For any given open neighborhood ~$\Omega_1$ ~of $\bdx_*$, 
there is a $\dl$ with $0<\dl<2$ such that 
$\overline{S_\dl(\bdx_*)}\subset\Omega_*\cap\Omega_1$  ~with
\begin{equation}\label{h12}
  \big\|J_\rkr(\bdx)^\dagger\big\|_2\,\big(\gamma\,\|\bdx-\bdy\|_2+\zeta\,
\|\bdf(\bdy)\|_2\big) ~~<~~ h ~~<~~ 1
\end{equation}
for all ~$\bdx,\,\bdy\,\in\,S_\dl(\bdx_*)$ ~and, there is a ~$\tau$ ~with 
$0<\tau<\frac{1}{2}\,\dl$ ~such that
\begin{equation}\label{al12}
\big\|J_\rkr(\bdz)^\dagger\big\|_2\,\|\bdf(\bdz)\|_2 ~~\le~~ 
\mbox{$\frac{1}{2}$}\,(1-h)\,\dl ~~<~~ \mbox{$\frac{1}{2}$}\,\dl ~~<~~ 1
\end{equation}
for all ~$\bdz\,\in\,S_\tau(\bdx_*)$.
Then, for every ~$\bdx_0\,\in\,S_\tau(\bdx_*)$, ~we have 
\begin{align*}
  \|\bdx_1-\bdx_*\| &~~\le~~ \|\bdx_1-\bdx_0\|_2+\|\bdx_0-\bdx_*\|_2 \\
&~~\le~~
\|J_\rkr(\bdx_0)^\dagger\|_2\,\|\bdf(\bdx_0)\|_2+\tau ~~<~~ \dl.
\end{align*}
Namely ~$\bdx_1\,\in\,S_\dl(\bdx_*)$.
Assume ~$\bdx_i\,\in\,S_\dl(\bdx_*)$ ~for all ~$i\,\in\,\{0,1,\ldots,k\}$.
Since
\begin{align*}
\bdx_j-\bdx_{j-1}+J_\rkr(\bdx_{j-1})^\dagger\, \bdf(\bdx_{j-1}) &~~=~~ \bdo
\tag{by (\ref{ni})} \nonumber \\
J_\rkr(\bdx_j)^\dagger\,J_\rkr(\bdx_j)\,J_\rkr(\bdx_j)^\dagger&~~=~~ 
J_\rkr(\bdx_j)^\dagger \tag{by (\ref{mpc})} \nonumber \\
J(\bdx_{j-1})\,J_\rkr(\bdx_{j-1})^\dagger &~~=~~
J_\rkr(\bdx_{j-1})\,J_\rkr(\bdx_{j-1})^\dagger
\tag{by (\ref{ara})}
\end{align*}
for all ~$j\,\in\,\{1,2,\ldots,k\}$, ~we have
\begin{eqnarray}
\lefteqn{\big\|\bdx_{j+1}-\bdx_j\big\|_2 
~~=~~  \big\|J_\rkr(\bdx_j)^\dagger\,\bdf(\bdx_j)\big\|_2} \nonumber  \\
&~=& \Big\|J_\rkr(\bdx_j)^\dagger\,\Big(\bdf(\bdx_j) -
J(\bdx_{j-1})\,
\big(\bdx_j-\bdx_{j-1}+J_\rkr(\bdx_{j-1})^\dagger\, \bdf(\bdx_{j-1})\big)
\Big)\Big\|_2  \nonumber \\
&~\le& \big\|J_\rkr(\bdx_j)^\dagger\,\big(\bdf(\bdx_j)-\bdf(\bdx_{j-1})
-J(\bdx_{j-1})\,(\bdx_j-\bdx_{j-1})\big)\big\|_2 \nonumber \\
& & + \big\|J_\rkr(\bdx_j)^\dagger\,\big(J_\rkr(\bdx_j)\,
J_\rkr(\bdx_j)^\dagger - 
J_\rkr(\bdx_{j-1})\,J_\rkr(\bdx_{j-1})^\dagger\big)\,\bdf(\bdx_{j-1})\big\|_2  
\nonumber \\
&~\le&  \|J_\rkr(\bdx_j)^\dagger\|_2\,\big(\gamma\,\|\bdx_j-\bdx_{j-1}\|_2 
+ \zeta\, \|\bdf(\bdx_{j-1})\|_2\big)
\,\|\bdx_j-\bdx_{j-1}\|_2  \label{6458}\\
&~\le& h\,\|\bdx_j-\bdx_{j-1}\|_2 \nonumber 
\end{eqnarray}
leading to
\[ \big\|\bdx_{j+1}-\bdx_j\big\|_2 ~~\le~~ h^j 
\|\bdx_1-\bdx_0\|_2 ~~\le~~ h^j\,\mbox{$\frac{1-h}{2}$}\,\dl
\]
for all ~$j\,\in\,\{1,\ldots,k\}$.
Thus
\begin{align*}
 \|\bdx_{k+1}-\bdx_0\|_2 &~~\le~~ \|\bdx_{k+1}-\bdx_k\|_2 + \cdots +
\|\bdx_1-\bdx_0\|_2 \\
&~~\le~~ \big(h^k+h^{k-1}+\cdots+1\big)\,
\|\bdx_1-\bdx_0\|_2 \\
&~~<~~ \mbox{$\frac{1}{1-h}$}\, \mbox{$\frac{1-h}{2}$}\,\dl
~~=~~\mbox{$\frac{1}{2}$}\,\dl
\end{align*}
and ~$\bdx_{k+1}\,\in\,S_\dl(\bdx_*)$ ~since
\[ \|\bdx_{k+1}-\bdx_*\|_2 ~~\le~~ \|\bdx_{k+1}-\bdx_0\|_2 + 
\|\bdx_0-\bdx_*\|_2  ~~<~~ \big(\mbox{$\frac{1}{2}+\frac{1}{2}$}\big)\,\dl
~~=~~ \dl,
\]
completing the induction so all iterates ~$\big\{\bdx_j\big\}_{j=0}^\infty$
of (\ref{ni}) are in ~$S_\dl(\bdx_*)$ ~from any initial iterate 
$\bdx_0\in S_\tau(\bdx_*)$.
Furthermore, the iterates ~$\big\{\bdx_j\big\}_{j=0}^\infty$ ~form a Cauchy 
sequence since, for any ~$k,\,j\,\ge\,0$, 
\begin{align}
\|\bdx_{k+j}-\bdx_k\|_2 &~~\le~~ \|\bdx_{k+j}-\bdx_{k+j-1}\|_2 +\cdots+
\|\bdx_{k+1}-\bdx_k\|_2 \nonumber \\
&~~\le~~ \big(h^{j-1}+\cdots+h+1\big)\,
\|\bdx_{k+1}-\bdx_k\|_2 \nonumber \\
&~~<~~ \mbox{$\frac{1}{1-h}$}\,\|\bdx_{k+1}-\bdx_k\|_2 \label{xkjk0} \\
&~~<~~ \mbox{$\frac{1}{1-h}$}\,h^k\cdot\mbox{$\frac{1-h}{2}$}\,\dl 
~~=~~ \mbox{$\frac{1}{2}$}\,\dl\,h^k
\label{xkjk}
\end{align}
can be as small as needed when ~$k$ ~is sufficient large.
Consequently, the sequence $\big\{\bdx_k\big\}_{k=0}^\infty$ generated
by the iteration (\ref{ni}) converges to a certain 
$\hat\bdx\,\in\,\overline{S_\dl(\bdx_*)}\subset\Omega_1$ ~at which 
$J_\rkr(\hat\bdx)^\dagger\,\bdf(\hat\bdx)\,=\,\bdo$ ~and thus, 
by Lemma~\ref{l:rz}, ~the limit ~$\hat\bdx$ ~is a semiregular zero in the same
branch of ~$\bdx_*$ ~with the identical dimension.

We now have 
\[ \|\bdx_k-\hat\bdx\|_2 ~~\le~~ \mbox{$\frac{1}{1-h}$}\,\|\bdx_k-\bdx_{k+1}
\|_2 ~~\le~~ \mbox{$\frac{1}{2}$}\,\dl\,h^k
\]
for all ~$k\,\ge\,1$ ~by setting 
$j\rightarrow\infty$ ~in (\ref{xkjk0}) and (\ref{xkjk}).
Substituting
\begin{align*} \|\bdf(\bdx_{j-1})\|_2 &~=~ 
\|\bdf(\bdx_{j-1})-\bdf(\hat\bdx)\|_2 ~\le~ 
\mu\,\|\bdx_{j-1}-\hat\bdx\|_2 \\
& ~\le~ \mu\,\big(\|\bdx_{j-1}-\bdx_j\|_2+\|\bdx_j-\bdx_{j+1}\|_2+
\|\bdx_{j+1}-\bdx_{j+2}\|_2+ \cdots\big) \\
& ~\le ~\frac{\mu}{1-h}\, \|\bdx_j-\bdx_{j-1}\|_2
\end{align*}
for a certain ~$\mu\,>\,0$, ~the inequality (\ref{6458}) yields
\[  \|\bdx_{j+1}-\bdx_j\|_2 ~~<~~ \bt\,\|\bdx_j-\bdx_{j-1}\|_2^2
\]
for all ~$j\,=\,1,2,\ldots$ where 
\[ \bt ~~=~~ \max\left\{1,\,\xi\, \left(\gamma + \frac{\zeta\,\mu}{1-h}\right) 
\right\}
\]
Let 
\[
\Omega_0 ~~=~~ \Big\{ \bdx\,\in\,S_\tau(\bdx_*) ~\Big|~
\big\|J_\rkr(\bdx)^\dagger \bdf(\bdx)\big\|_2 < \mbox{$\frac{1}{\bt}$}\,h 
\Big\}
\]
Then, for every initial iterate ~$\bdx_0\,\in\,\Omega_0$, ~we have
\[  \|\bdx_1-\bdx_0\|_2 ~~=~~\big\|J_\rkr(\bdx_0)^\dagger \bdf(\bdx_0)\big\|_2 
~~<~~ \mbox{$\frac{1}{\bt}$}\,h 
~~\le~~ h.
\]
Consequently, for all $j=1,2,\ldots$, ~we have
\begin{align*}
 \|\bdx_{j+1}-\bdx_j\|_2 &~~<~~ \bt^{2^j-1} \|\bdx_1-\bdx_0\|_2^{2^j}\\
&~~=~~ \big(\bt\,\|\bdx_1-\bdx_0\|_2\big)^{2^j-1}\,\|\bdx_1-\bdx_0\|_2
~~<~~ h^{2^j}
\end{align*}
Thus
\begin{align}  \|\bdx_k-\hat\bdx\|_2 &~~=~~ 
\lim_{j\rightarrow\infty} \|\bdx_{k+j}-\bdx_k\|_2 \nonumber \\
& ~~\le~~ 
\lim_{j\rightarrow\infty} \big(\|\bdx_{k+j}-\bdx_{k+j-1}\|_2 + \cdots +
\|\bdx_{k+1}-\bdx_k\|_2\big) \nonumber  \\
&~~\le~~ \mbox{$\frac{1}{1-h}$}\,\|\bdx_{k+1}-\bdx_k\|_2
~~\le~~ 
\mbox{$\frac{1}{1-h}$}\,h^{2^k}
\label{qci}
\end{align}
with ~$h\,<\,1$. 
Consequently the convergence to ~$\hat\bdx$ ~is at quadratic rate. 
~\qed

Theorem~\ref{t:mt} can serve as the universal convergence theorem 
of Newton's iteration including the conventional version (\ref{ni0}) as
a special case.
The sufficient conditions for Theorem~\ref{t:mt} consist of 
smoothness of $\bdf$, semiregularity of $\bdx_*$ and the initial iterate
$\bdx_0$ being near $\bdx_*$.
These assumptions are minimal and indispensable for the convergence of 
the rank-$r$ Newton's iteration (\ref{ni}). 
The iteration would need to be adapted if the mapping $\bdf$ is not 
smooth.
Without the semiregularity of $\bdx_*$, the limit $\hat\bdx$ as a 
stationary point is seldom a zero from our experiment.
The non-global convergence is an accepted imperfection of Newton's iteration 
rather than a drawback.
Theoretical criteria on the initial iterate beyond sufficient nearness to
a solution are of little practical meaning as trial-and-error would be 
easier than verifying those conditions.

The common definition of quadratic convergence of a sequence
$\{\bdx_k\}$ to its limit $\hat\bdx$ would require that there is a 
constant $\la$ such that $\|\bdx_k-\hat\bdx\|_2 \,\le\,\la\,\|\bdx_{k-1}-
\hat\bdx\|_2^2$ for $k$ larger than a certain $k_0$ and imply 
$\|\bdx_k-\hat\bdx\|_2 \,\le\, \la\,\|\bdx_{k_0}-\hat\bdx\|_2^{2^{k-k_0}}$.
The inequality (\ref{qci}) ensures essentially the same error bound for 
the iterate $\bdx_k$ toward $\hat\bdx$ and can be used as an alternative
definition of quadratic convergence.

\begin{remark}[Convergence near an ultrasingular zero]\em
If $\hat\bdx$ is an ultrasingular zero of $\bdf$ where 
$r\,=\,\rank{J(\hat\bdx)}$, the iteration (\ref{ni}) still locally converges 
to a certain stationary point $\check\bdx$ at which 
$J_\rkr(\check\bdx)^\dagger\,\bdf(\check\bdx)\,=\,\bdo$ but 
$\check\bdx$ is not necessarily a zero of $\bdf$.
Such an $\check\bdx$ satisfies the necessary condition for 
$\|\bdf(\bdx)\|_2$ to reach a local minimum if $\rank{J(\check\bdx)}\,=\,r$. 
From our computational experiment, a stationary point to which the 
iteration (\ref{ni}) converges is rarely a zero of $\bdf$ 
when the initial iterate is near an ultrasingular zero.
\end{remark}

\begin{remark}[On the projection rank $r$]\em
Application of the iteration (\ref{ni}) requires identification of the 
rank $r$ of the Jacobian at a zero.
There are various approaches for determining $r$ such as analytical
methods on the application model (c.f. \S\ref{s:gcd} and \S\ref{s:eig}),
numerical matrix rank-revealing \cite{golub-vanloan4,lee-li-zeng,li-zeng-03}, 
and even trial-and-error.
For any positive error tolerance 
$\theta\,<\,\big\|J(\bdx_*)^\dagger\|_2^{-1}$, the numerical rank
$\ranka{\theta}{J(\bdx_0)}$ within $\theta$ is identical to 
$r\,=\,\rank{J(\bdx_*)}$
if $\bdx_0$ is sufficiently close to $\bdx_*$ \cite{Zeng2019}.
Furthermore, the projection rank needs to be identified or computed only 
once for a solution branch and the same rank can be used repeatedly 
in calculating other witness points in the same branch.
\end{remark}

\section{Convergence theorem on perturbed equations}\label{s:cp}

Practical applications in scientific computing often involve equations that 
are given through empirical data with limited accuracy.
While the underlying exact equation can have solution sets of positive
dimensions, the perturbed equation may not. 
Such applications can be modeled as solving an equation
\begin{equation}\label{fxy0}
\bdf(\bdx,\,\bdy) ~~=~~ \bdo ~~~~\mbox{for}~~~~ \bdx\,\in\,\Omega\subset\,\C^m
~~\mbox{or}~~ \R^m
\end{equation}
at a particular parameter value $\bdy\in\Sigma\subset\C^n$ \,or \,$\R^n$ 
representing the data.
The equation (\ref{fxy0}) may have a nonisolated solution set only at a 
particular isolated parameter value $\bdy\,=\,\bdy_*$.
A typical example is given as Example~\ref{e:cyc} later in \S\ref{s:nag} 
where the 1-dimensional solutions in $\bdx$ exist for an equation 
$\bdf(\bdx,t)\,=\,\bdo$ only when $t\,=\,1$ exactly. 
The question becomes: {\em
Assuming the equation {\em (\ref{fxy0})} has a semiregular solution 
$\bdx\,=\,\bdx_*$ at an underlying parameter $\bdy\,=\,\bdy_*$ but 
$\bdy_*$ is known only through empirical data in 
$\tilde\bdy\,\approx\,\bdy_*$, does the iteration
\begin{equation}\label{niy}
\bdx_{k+1} ~~=~~ \bdx_k - \bdf_\bdx(\bdx_k,\,\tilde\bdy)_\rkr^\dagger\,
\bdf(\bdx_k,\,\tilde\bdy), ~~~~k\,=\,0,1,\ldots
\end{equation}
converge and, if so, does the limit $\tilde\bdx$ approximate a zero 
$\bdx\,=\,\hat\bdx$ of the mapping $\bdx\,\mapsto\,\bdf(\bdx,\bdy_*)$ 
with an accuracy $\|\tilde\bdx-\hat\bdx\|_2 \,=\, O(\|\tilde\bdy-\bdy_*\|_2)$
in the same order of the data?}
The following theorem attempts to answer that question.

\begin{theorem}[Convergence Theorem on Perturbed Equations]\label{t:mt2}
\!\!Let a smooth mapping $(\bdx,\,\bdy)\,\mapsto\,\bdf(\bdx,\,\bdy)$ 
be defined in an open domain.
Assume $\bdx_*$ is a semiregular zero of
$\bdx\,\mapsto\,\bdf(\bdx,\bdy_*)$ at a fixed 
$\bdy_*$ with $\rank{\bdf_\bdx(\bdx_*,\bdy_*)}\,=\,r\,>\,0$ and
$\|\bdf_\bdy(\bdx_*,\bdy_*)\|_2\,>\,0$.
Then there exist a neighborhood $\Omega_*\times\Sigma_*$ of 
$(\bdx_*,\bdy_*)$, a neighborhood $\Omega_0$ of $\bdx_*$ and a constant
$h$ with $0<h<1$ such that, 
at every fixed $\tilde\bdy\in\Sigma_*$ serving as empirical data for
$\bdy_*$ and from any initial iterate $\bdx_0\in\Omega_0$, 
the iteration {\em (\ref{niy})} linearly converges to a stationary point 
$\tilde\bdx\in\Omega_*$ at which 
$\bdf_\bdx(\tilde\bdx,\tilde\bdy)_\rkr^\dagger\,\bdf(\tilde\bdx,\tilde\bdy)
\,=\,\bdo$ with an error bound
\begin{align}\label{eby}
\|\tilde\bdx-\hat\bdx\|_2 &~~\le~~  \mbox{$\frac{8}{1-h}$}\,
\big\|\bdf_\bdx(\bdx_*,\bdy_*)^\dagger\big\|_2\,
\big\|\bdf_\bdy(\bdx_*,\bdy_*)\big\|_2\,\|\tilde\bdy-\bdy_*\|_2 
+O\big(\|\tilde\bdy-\bdy_*\|_2^2\big) 
\end{align}
toward a semiregular zero $\hat\bdx$ of $\bdx\,\mapsto\,\bdf(\bdx,\bdy_*)$
in the same branch of $\bdx_*$.
\end{theorem}

{\em Proof.}
Following almost the same proof of Lemma~\ref{l:p} using the smoothness of
the mapping $(\bdx,\,\bdy)\,\mapsto\,\bdf(\bdx,\,\bdy)$, there is an open 
bounded convex neighborhood $\Omega_1\times\Sigma_1$ of $(\bdx_*,\bdy_*)$ 
along with constants $\zeta,\,\gamma,\,\eta\,>\,0$ such that
\begin{align*}  
\big\|\bdf_\bdx(\hat\bdx,\hat\bdy)_\rkr\, \bdf_\bdx(\hat\bdx,\hat\bdy
)_\rkr^\dagger - \bdf_\bdx(\check\bdx,\hat\bdy)_\rkr\, 
\bdf_\bdx(\check\bdx,\hat\bdy)_\rkr^\dagger\big\|_2 &~~\le~~
\zeta\,\big\|\hat\bdx-\check\bdx\big\|_2 \\
\big\|\bdf(\hat\bdx,\hat\bdy)-\bdf(\check\bdx,\hat\bdy)-\bdf_\bdx(\check\bdx,
\hat\bdy)\,
(\hat\bdx-\check\bdx)\big\|_2 &~~\le~~ \gamma\,
\big\|\hat\bdx-\check\bdx\big\|_2^2  \\
\big\|\bdf_\bdx(\hat\bdx,\hat\bdy)_\rkr^\dagger - 
\bdf_\bdx(\hat\bdx,\check\bdy)_\rkr^\dagger\big\|_2 &~~\le~~ \eta\,
\big\|\hat\bdy-\check\bdy\big\|_2 \\
\big\|\bdf_\bdx(\hat\bdx,\hat\bdy)\big\|_2 ~~<~~ 2\,
\big\|\bdf_\bdx(\bdx_*,\bdy_*)\big\|_2 ~~~\mbox{and}~~~
\big\|\bdf_\bdy(\hat\bdx,\hat\bdy)\big\|_2 &~~<~~ 2\,
\big\|\bdf_\bdy(\bdx_*,\bdy_*)\big\|_2 
\end{align*}
for all $\hat\bdx,\,\check\bdx\in\Omega_1$ and 
$\hat\bdy,\,\check\bdy\in\Sigma_1$ with
\[
\max_{(\hat\bdx,\hat\bdy)\in\overline{\Omega}_1\times\overline{\Sigma}_1} 
\big\|\bdf_\bdx(\hat\bdx,\hat\bdy)_\rkr^\dagger\big\|_2 ~~\le~~
2\,\big\|\bdf_\bdx(\bdx_*,\bdy_*)^\dagger\big\|_2
\]
Let 
$S_\eps(\bdx_*)\,:=\,\big\{\bdx\in\Omega_1~\big|~\|\bdx-\bdx_*\|_2<\eps\big\}$
for any $\eps>0$ and $h$ be any fixed constant with $0<h<1$.
There are constants $\dl,\,\tau,\,\tau'>0$ with 
$6\,\tau'\,<\,2\,\tau\,<\,\dl$,
$S_\tau(\bdx_*)\subset S_\dl(\bdx_*)\subset\Omega_1$ and an open neighborhood
$\Sigma_0\subset\Sigma_1$ of $\bdy_*$ such that
\begin{align*}
  \big\|\bdf_\bdx(\hat\bdx,\hat\bdy)_\rkr^\dagger\big\|_2\,
\big(\gamma\,\|\hat\bdx
-\check\bdx\|_2 + \zeta\,\|\bdf(\hat\bdx,\hat\bdy)\|_2\big) &~~<~~ h \\
\big\|\bdf_\bdx(\bdz,\hat\bdy)_\rkr^\dagger\big\|_2\,\|\bdf(\bdz,\hat\bdy)\|_2
&~~\le~~ \mbox{$\frac{1}{2}$}\,(1-h)\,\dl ~<~ \mbox{$\frac{1}{2}$}\,\dl \\
\big\|\bdf_\bdx(\tilde\bdz,\hat\bdy)_\rkr^\dagger\big\|_2\,\|
\bdf(\tilde\bdz,\hat\bdy)\|_2
&~~\le~~ \mbox{$\frac{1}{3}$}\,(1-h)\,\tau ~<~ \mbox{$\frac{1}{3}$}\,\tau
\end{align*}
for all $\hat\bdx,\,\check\bdx\in S_\dl(\bdx_*)$, $\bdz\in S_\tau(\bdx_*)$, 
$\tilde\bdz\in S_{\tau'}(\bdx_*)$ and $\hat\bdy\in\Sigma_0$.
Using the same argument in the proof of Theorem~\ref{t:mt}, the sequence
$\{\bdx_k\}$ generated by the iteration (\ref{niy}) starting from any 
$\bdx_0\in S_{\tau'}(\bdx_*)$ satisfies
$\|\bdx_{k+1}-\bdx_k\|_2\,<\,h\,\|\bdx_k-\bdx_{k-1}\|_2$ for all 
$k\,\ge\,1$ and is a Cauchy sequence 
staying in $S_\tau(\bdx_*)$ for every fixed $\tilde\bdy\in\Sigma_0$ 
and converges to an 
$\tilde\bdx\in \overline{S_{\frac{2}{3}\tau}(\bdx_*)}\subset S_\tau(\bdx_*)$ 
that depends on the choices of 
$\bdx_0$ and $\tilde\bdy$.
Resetting $\bdx_0\,=\,\tilde\bdx\in S_\tau(\bdx_*)$, 
Theorem~\ref{t:mt} ensures the iteration (\ref{niy}) with 
$\bdy\,=\,\bdy_*\in\Sigma_0$ stays in $S_\dl(\bdx_*)$ and converges to a 
certain semiregular zero $\hat\bdx$ of the mapping
$\bdx\,\mapsto\,\bdf(\bdx,\bdy_*)$.
We can assume 
\[  \big\|\hat\bdy-\check\bdy\big\|_2 ~~<~~ 
\frac{1-h}{4\,\eta\,\|\bdf_\bdx(\bdx_*,\bdy_*)\|_2}
\]
for all $\hat\bdy,\,\check\bdy\in\Sigma_0$ by shrinking $\Sigma_0$ if
necessary.
Subtracting both sides of
\begin{align*}
\tilde\bdx~ &~~=~~ \tilde\bdx - \bdf_\bdx(\tilde\bdx,\tilde\bdy)_\rkr^\dagger\,
\bdf(\tilde\bdx,\tilde\bdy) \\
\bdx_1 &~~=~~ \tilde\bdx - \bdf_\bdx(\tilde\bdx,\bdy_*)_\rkr^\dagger\, 
\bdf(\tilde\bdx,\bdy_*)
\end{align*}
yields
\begin{align*} 
\|\tilde\bdx-\bdx_1\|_2 &~~=~~ \big\|
\bdf_\bdx(\tilde\bdx,\tilde\bdy)_\rkr^\dagger\,\bdf(\tilde\bdx,\tilde\bdy) -
\bdf_\bdx(\tilde\bdx,\bdy_*)_\rkr^\dagger\, \bdf(\tilde\bdx,\bdy_*)\big\|_2 \\
&~~\le~~ \big\|\bdf_\bdx(\tilde\bdx,\tilde\bdy)_\rkr^\dagger - 
\bdf_\bdx(\tilde\bdx,\bdy_*)_\rkr^\dagger\big\|_2\,
\|\bdf(\tilde\bdx,\tilde\bdy)-\bdf(\hat\bdx,\bdy_*) \|_2 \\
&~~~~~~ + \big\|\bdf_\bdx(\tilde\bdx,\bdy_*)_\rkr^\dagger\big\|_2\,
\|\bdf(\tilde\bdx,\tilde\bdy)-\bdf(\tilde\bdx,\bdy_*)\|_2.
\end{align*}
From 
\begin{align*} 
\lefteqn{\big\|\bdf_\bdx(\tilde\bdx,\bdy_*)_\rkr^\dagger\big\|_2\,
\|\bdf(\tilde\bdx,\tilde\bdy)-\bdf(\tilde\bdx,\bdy_*)\|_2} \\
&~~\le~~ 
4\,\big\|\bdf_\bdx(\bdx_*,\bdy_*)^\dagger\big\|_2\, 
\big\|\bdf_\bdy(\bdx_*,\bdy_*)\big\|_2\,\|\tilde\bdy-\bdy_*\|_2
\end{align*}
and, since $\bdf(\hat\bdx,\bdy_*)\,=\,\bdo$,
\begin{align*}
\lefteqn{\big\|\bdf_\bdx(\tilde\bdx,\tilde\bdy)_\rkr^\dagger -
\bdf_\bdx(\tilde\bdx,\bdy_*)_\rkr^\dagger\big\|_2\,\|\bdf(\tilde\bdx,\tilde\bdy)
-\bdf(\hat\bdx,\bdy_*)\|_2} \\
&~~<~~ \eta\,\big\|\tilde\bdy-\bdy_*\big\|_2 \,\big(
2\,\|\bdf_\bdx(\bdx_*,\bdy_*)\|_2\,\|\tilde\bdx-\hat\bdx\|_2 + 
O(\|\tilde\bdy-\bdy_*\|_2)
\big) \\
&~~\le~~ \mbox{$\frac{1-h}{2}$} \,\|\hat\bdx-\tilde\bdx\|+
O\big(\|\tilde\bdy-\bdy_*\|_2^2\big),
\end{align*}
we have 
\[ \|\tilde\bdx-\bdx_1\|_2 ~~\le~~ 
4\,\big\|\bdf_\bdx(\bdx_*,\bdy_*)^\dagger\big\|_2\, 
\big\|\bdf_\bdy(\bdx_*,\bdy_*)\big\|_2\,\|\tilde\bdy-\bdy_*\|_2
+ \mbox{$\frac{1-h}{2}$} \,\|\hat\bdx-\tilde\bdx\|+
O\big(\|\tilde\bdy-\bdy_*\|_2^2\big).
\]
As a result, 
\begin{align*}
\lefteqn{\|\hat\bdx-\tilde\bdx\|_2 ~~=~~ 
\lim_{k\rightarrow\infty}  \|\bdx_k-\tilde\bdx\|_2}  \\
 &~~\le~~ \lim_{k\rightarrow\infty}\big(\|\bdx_k-\bdx_{k-1}\|_2+\cdots+
\|\bdx_1-\tilde\bdx\|_2\big) \\
& ~~\le~~ \lim_{k\rightarrow\infty}(h^k+h^{k-1}+\cdots+1)\,
\|\bdx_1-\tilde\bdx\|_2 \\
&~~\le~~ 
\mbox{$\frac{4}{1-h}$}\,\big\|\bdf_\bdx(\bdx_*,\bdy_*)^\dagger\big\|_2\, 
\big\|\bdf_\bdy(\bdx_*,\bdy_*)\big\|_2\,\|\tilde\bdy-\bdy_*\|_2
+ \mbox{$\frac{1}{2}$} \,\|\hat\bdx-\tilde\bdx\|+
O\big(\|\tilde\bdy-\bdy_*\|_2^2\big),
\end{align*}
leading to (\ref{eby}). 
The theorem is proved by setting 
$\Omega_*\times\Sigma_*\,=\,S_\dl(\bdx_*)\times\Sigma_0$ and 
$\Omega_0\,=\,S_{\tau'}(\bdx_*)$.
~\qed

Theorem~\ref{t:mt2} leads to some intriguing implications.
At the exact data parameter value $\bdy\,=\,\bdy_*$, solving the equation 
$\bdf(\bdx,\bdy_*)\,=\,\bdo$ for a nonisolated solution in $\bdx$ 
can be an ill-posed problem as the structure of the solution can be 
infinitely sensitive to perturbations on the parameter $\bdy$.
However, Theorem~\ref{t:mt2} implies that solving the stationary equation 
\[  \bdf_\bdx(\bdx,\bdy_*)_\rk{r}^\dagger\,\bdf(\bdx,\bdy_*) ~~=~~ \bdo
~~~~\mbox{for}~~~ \bdx\,\in\,\Omega
\]
in a neighborhood of a semiregular zero $\bdx_*$ is a well-posed problem in
the sense that the solution is Lipschitz continuous with 
respect to the perturbations on the parameter $\bdy$ from $\bdy_*$
with $r\,=\,\rank{\bdf_\bdx(\bdx_*,\bdy_*)}$.
When the parameter $\bdy_*$ is only known through empirical data 
$\tilde\bdy$, the mapping $\bdx\,\mapsto\,\bdf(\bdx,\tilde\bdy)$ may
not have a positive-dimensional zero of its own near $\bdx_*$ and, in 
some cases, does not have a zero at all.
As a result, solving the equation $\bdf(\bdx,\tilde\bdy)\,=\,\bdo$ 
in exact sense for $\bdx$ is futile even if we extend the machine precision. 
It may be a pleasant surprise that a zero $\tilde\bdx$
of the mapping $\bdx\,\mapsto\,\bdf_\bdx(\bdx,\tilde\bdy)_\rk{r}^\dagger\,
\bdf(\bdx,\tilde\bdy)$ exists near $\bdx_*$. 
Furthermore that $\tilde\bdx$ approximates a desired zero
$\hat\bdx$ of the underlying mapping $\bdx\,\mapsto\,\bdf(\bdx,\bdy_*)$ 
with an accuracy in the same order as accuracy of the data.
More importantly in practice, this numerical zero $\tilde\bdx$ is 
attainable as the rank-$r$ Newton's iteration (\ref{niy}) locally converges 
to it.

The iteration (\ref{niy}) is more than an algorithm as it is, at the same time,
a natural regularization of a hypersensitive zero-finding problem. 
Since $\bdx_*$ is a nonisolated zero of the mapping 
$\bdx\,\mapsto\,\bdf(\bdx,\bdy_*)$, the iteration does not necessarily 
converge to $\bdx_*$ but to some $\hat\bdx$ in the same branch 
of the solution set.

As a special case, Theorem~\ref{t:mt2} reaffirms the local convergence 
of the Gauss-Newton iteration in solving an overdetermined system 
$\bdf(\bdx)\,=\,\bdb$ for its least squares solution%
\cite[Theorem 10.2.1, p.222]{dennis-schnabel}\cite{ZengAIF}
and extends the result with an error bound on the solution 
with respect to the perturbation on the right-hand side $\bdb$.

\begin{corollary}[Convergence of Gauss-Newton Iteration]
Let $\bdx\,\mapsto\,\bdf(\bdx)$ be a smooth mapping from an open domain 
$\Omega\subset\F^n$ to $\F^m$ where $\F\,=\,\R$ or $\C$ with $n\,<\,m$.
Assume $\bdx_*\,\in\,\Omega$ is a least squares solution of the 
overdetermined system $\bdf(\bdx_*)\,=\,\bdb$ with 
$\nullity{\bdf_\bdx(\bdx_*)}\,=\,0$.
Then there are open neighborhoods $\Omega_*$ and $\Sigma_*$
of $\bdx_*$ and $\bdf(\bdx_*)$ respectively such that the Gauss-Newton
iteration
\[ \bdx_{k+1} ~~=~~ \bdx_k - \bdf_\bdx(\bdx_k)^\dagger(\bdf(\bdx_k)-\bdb)
~~~~\mbox{for}~~~ k\,=\,0,1, \cdots
\]
linearly converges to $\bdx_*$ from any initial iterate 
$\bdx_0\,\in\,\Omega_*$ provided that $\bdb\,\in\,\Sigma_*$ and the 
convergence is quadratic if $\bdb\,=\,\bdf(\bdx_*)$.
Furthermore, the same iteration 
\[ \bdx_{k+1} ~~=~~ \bdx_k - \bdf_\bdx(\bdx_k)^\dagger(\bdf(\bdx_k)-\tilde\bdb)
~~~~\mbox{for}~~~ k\,=\,0,1, \cdots
\]
on the perturbed system $\bdf(\bdx)\,=\,\tilde\bdb\,\in\,\Sigma_*$ 
from $\bdx_0\,\in\,\Omega_*$
converges to a critical point $\tilde\bdx \,\in\,\Omega$ at which 
$\bdf_\bdx(\tilde\bdx)^\h\,\big(\bdf(\tilde\bdx)-\tilde\bdb\big)\,=\,\bdo$
with an error bound
\[ \|\tilde\bdx-\bdx_*\|_2 ~~\le ~~ 
c\,\big\|\bdf_\bdx(\bdx_*)^\dagger\big\|_2\,\big\|\tilde\bdb-\bdb\big\|_2 
+ O\big(\|\tilde\bdb-\bdb\|_2^2\big)
\]
where $c$ is a moderate constant depending on $\tilde\bdx$ and 
$\tilde\bdb$.
\end{corollary}

{\em Proof.} The proof is a straightforward verification by applying 
Theorem~\ref{t:mt2} on the mapping $(\bdx,\,\bdy)\,\mapsto\,\bdf(\bdx)-\bdy$
at $(\bdx_*,\,\bdy_*)$ where $\bdy_*\,=\,\bdf(\bdx_*)$ with 
$\tilde\bdy\,=\,\bdb$ and $\hat\bdx\,=\,\bdx_*$. \qed

\begin{remark}[Condition number of a nonisolated zero]\em
\,\,Another implication of Theorem~\ref{t:mt2} lies in the sensitivity measure
derived from the error estimate (\ref{eby}), from which we can naturally 
define a {\em condition number} of a zero $\bdx_*$ of the mapping
$\bdx\,\mapsto\,\bdf(\bdx,\bdy_*)$ with respect to the parameter value 
$\bdy\,=\,\bdy_*$ as
\begin{equation}\label{cdn}
\kappa_{\bdf,\bdx}(\bdx_*,\bdy_*) ~~:=~~ \left\{\begin{array}{ccl}
\big\|\bdf_\bdx(\bdx_*,\bdy_*)^\dagger\big\|_2\,\|\bdf_\bdy(\bdx_*,\bdy_*)\|_2
& & \mbox{if ~$\bdx_*$ ~is semiregular} \\
\infty & & \mbox{otherwise.} \end{array}\right.
\end{equation}
The condition number of a semiregular zero is finite and, if it is not large, 
a small perturbation in the data parameter $\bdy$ results in accurate 
approximation of the zero with an error estimate (\ref{eby}). 
On the other hand, errors of an ultrasingular zero has no known bound and the
condition number is thus infinity.
\end{remark}

\section{A geometric interpretation}\label{s:geo}

We consider a special case first: Solving a consistent linear system 
$A\,\bdx\,=\,\bdb$ with a rank $r$ matrix $A\,\in\,\C^{m\times n}$
using the rank-$r$ Newton's iteration (\ref{ni}). 
Since the system is consistent, namely $\bdb\,\in\,\cR(A)$, the 
solution set is an $(n-r)$-dimensional affine subspace 
\begin{align*}  A^\dagger\,\bdb+\cK(A) &~~:=~~ 
\big\{A^\dagger\,\bdb+\bdy ~\big|~ \bdy\,\in\,\cK(A)\big\} \\
&~~=~~ 
\big\{A^\dagger\,\bdb+N\,\bdz ~\big|~ \bdz\,\in\,\C^{n-r}\big\} 
\end{align*}
in which every particular solution is semiregular as defined in 
Definition~\ref{d:rz} where columns of 
$N\,\in\,C^{n\times(n-r)}$ form an orthonormal basis for $\cK(A)$.
From any initial iterate $\bdx_0\,\in\,\C^n$, the nearest point in 
the solution set
$A^\dagger\,\bdb+\cK(A)$ is $A^\dagger\,\bdb+N\,\bdz_0$ where 
$\bdz\,=\,\bdz_0$ is the least squares solution to the linear system
\[  A^\dagger\,\bdb+N\,\bdz ~~=~~ \bdx_0.
\]
Namely $\bdz_0 \,=\, N^\h\,(\bdx_0-A^\dagger\,\bdb) \,=\, N^\h\,\bdx_0$
since $N^\h\,A^\dagger\,=\,O$, and the nearest solution
\[  A^\dagger\,\bdb+N\,\bdz_0 ~~=~~ A^\dagger\,\bdb+N\,N^\h\,\bdx_0 
~~=~~  A^\dagger\,\bdb+(I-A^\dagger\,A)\,\bdx_0 
\]
since $(I-A^\dagger\,A)\,=\,N\,N^\h$ are the same 
orthogonal projection onto $\cK(A)$.
On the other hand, let the mapping $\bdf\,:\,\C^n\,\rightarrow\,\C^m$ be 
defined as $\bdf(\bdx)\,=\,A\,\bdx-\bdb$ with the Jacobian
$J(\bdx)\,\equiv\,A\,\equiv\,J_\rkr(\bdx)$.
From $\bdx_0\,\in\,\C^n$, the rank-$r$ Newton's iteration (\ref{ni}) 
requires only one step
\[  \bdx_1 ~~=~~ \bdx_0 - A^\dagger\,(A\,\bdx_0-\bdb)
~~=~~ A^\dagger\,\bdb + (I-A^\dagger\,A)\,\bdx_0
\]
In other words, the rank-$r$ Newton's iteration converges to the nearest
solution on the $(n-r)$-dimensional solution set from the initial iterate.

We now consider a general nonlinear mapping $\bdf\,:\,\Omega\subset\C^m\,
\rightarrow\,\C^n$ with the Jacobian $J(\bdx)$ at any $\bdx\in\Omega$.
At an iterate $\bdx_j$ near a semiregular $(n-r)$-dimensional
solution $\hat\bdx$ of the equation $\bdf(\bdx)\,=\,\bdo$, 
we have 
\begin{align*}   \bdf(\bdx) &~~=~~ \bdf(\bdx_j)+J(\bdx_j)\,(\bdx-\bdx_j) + 
O(\|\bdx-\bdx_j\|_2^2) \\
&~~\approx~~ \bdf(\bdx_j)+J_\rkr(\bdx_j)\,(\bdx-\bdx_j).
\end{align*}
On the other hand, from 
$\bdf(\bdx_j)\,=\, J(\hat\bdx)\,(\bdx_j-\hat\bdx) 
+ O\big(\|\bdx_j-\hat\bdx\|_2^2\big)$
we have 
\[ \big(I-J(\hat\bdx)\,J(\hat\bdx)^\dagger\big)\,
\big(\bdf(\bdx_j)-\bdf(\hat\bdx)\big) ~~=~~ O\big(\|\bdx_j-\hat\bdx\|_2^2\big)
\]
and
\begin{align*}
\bdf(\bdx_j) &~~=~~ 
J_\rkr(\bdx_j)\,J_\rkr(\bdx_j)^\dagger\,\bdf(\bdx_j) +
\big(I-J_\rkr(\bdx_j)\,J_\rkr(\bdx_j)^\dagger\big)\,\bdf(\bdx_j) \\
&~~=~~
J_\rkr(\bdx_j)\,J_\rkr(\bdx_j)^\dagger\,\bdf(\bdx_j)  +
\big(I-J(\hat\bdx)\,J(\hat\bdx)^\dagger\big)\, 
\big(\bdf(\bdx_j)-\bdf(\hat\bdx)\big) \\
&~~~~~~~~~~~~~
+\big(J(\hat\bdx)\,J(\hat\bdx)^\dagger-J_\rkr(\bdx_j)\,J_\rkr(\bdx_j)^\dagger
\big)\, \big(\bdf(\bdx_j) - \bdf(\hat\bdx)\big)\\
&~~=~~ 
J_\rkr(\bdx_j)\,J_\rkr(\bdx_j)^\dagger\,\bdf(\bdx_j) + 
O\big(\|\bdx_j-\hat\bdx\|_2^2\big)
\end{align*}
The basic principle of Newton's iteration is to designate the
next iterate $\bdx\,=\,\bdx_{j+1}$ as a numerical solution of the linear
system
\begin{equation}\label{nteq}   
\bdf(\bdx_j)+J(\bdx_j)\,(\bdx-\bdx_j) ~~=~~ \bdo.
\end{equation}
Since $J(\hat\bdx)$ is rank-deficient at the nonisolated zero 
$\hat\bdx$ and $\bdx_j$ is near $\hat\bdx$, the coefficient matrix
$J(\bdx_j)$ of the system (\ref{nteq}) is expected to be highly 
ill-conditioned for being nearly rank-deficient.
Consequently, solving the system (\ref{nteq}) as it is may not be reliable.
As elaborated in \cite[\S 8]{Zeng2019}, an alternatives is to solve
the nearby linear system
\begin{equation}\label{nteq1}
 J_\rkr(\bdx_j)\,J_\rkr(\bdx_j)^\dagger\,\bdf(\bdx_j)  
+J_\rkr(\bdx_j)\,(\bdx-\bdx_j) ~~=~~ \bdo
\end{equation}
and pick a proper vector $\bdx\,=\,\bdx_{j+1}$ from the solution
in the affine Grassmannian.
The linear system (\ref{nteq1}) is well-conditioned if the actual sensitivity 
measure
$\|J(\hat\bdx)\|_2\,\|J(\hat\bdx)^\dagger\|_2\,\approx\,
\|J_\rkr(\bdx_j)\|_2\,\|J_\rkr(\bdx_j)^\dagger\|_2$ is 
not large and (\ref{nteq1}) is an approximation to the system (\ref{nteq}).
The iterate $\bdx\,=\,\bdx_{j+1}$ from (\ref{ni}) is the minimum norm 
solution of the approximate system (\ref{nteq1}) and 
\[  \bdx_{j+1}-\bdx_j ~~\in~~ \cK(J_\rkr(\bdx_j))^\perp 
~~\approx~~ \cK(J(\hat\bdx))^\perp ~~=~~ \cR(\phi_\bdz(\hat\bdz))^\perp
\]
where $\bdz\,\mapsto\,\phi(\bdz)$ is the parametrization of the 
$(n-r)$-dimensional zero of $\bdf$ with $\phi(\hat\bdz)\,=\,\hat\bdx$,
implying the geometric interpretation of the rank-$r$ Newton's 
iteration:
\begin{quote}
From an initial iterate, the rank-$r$ Newton's iteration (\ref{ni}) aims at 
the nearest point on the solution branch following a normal line toward the 
solution set.
\end{quote}
How accurate the iteration hitting the nearest point depending on several
factors including how far the initial iterate is.
Similar geometric interpretations are also observed in \cite{Chu83,Tenabe}
in the cases where the Jacobians are surjective.

\begin{figure}[ht]
\begin{center}
\epsfig{figure=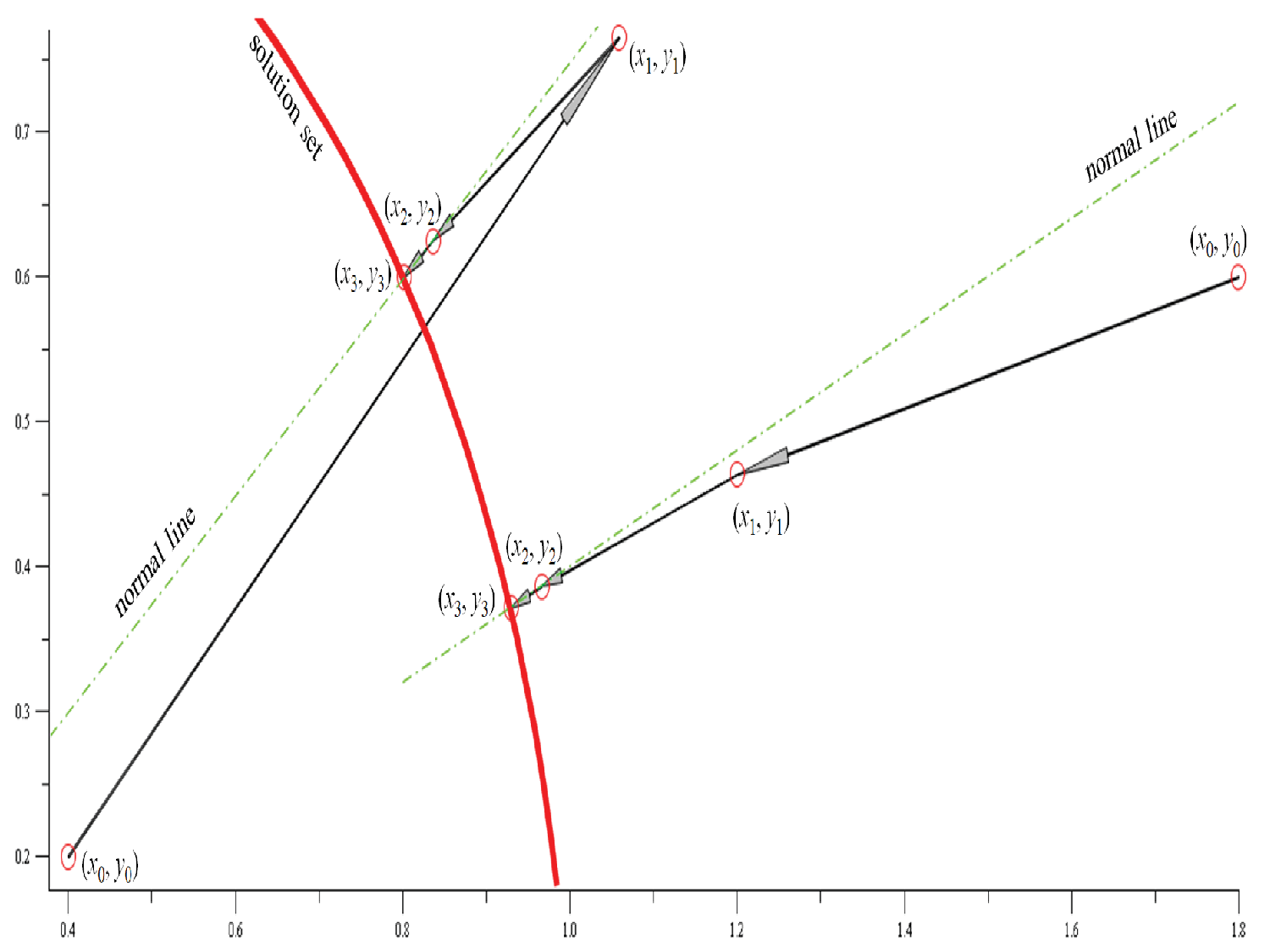,width=5in,height=2.2in}
\end{center}
\caption{\small Each sequence of iteration (\ref{ni}) asymptotically follows 
a normal line toward the solution set. Illustration is plotted using
actual data in Example~\ref{e:dir} from two initial iterates.
}\label{f:dir}
\end{figure}

\begin{example}[Normal line direction of convergence]\label{e:dir}\em
Consider $\bdf\,:\,\R^2\,\rightarrow\,\R^2$ as follows
\begin{equation}\label{f2}
\bdf(x,y) ~~=~~ \left(\begin{array}{c}
x^3+x\,y^2-x+2\,x^2+2\,y^2-2 \\ x^2\,y+y^3-y-3\,x^2-3\,y^2+3
\end{array}\right)
\end{equation}
whose zeros consist of a semiregular 1-dimensional unit circle $x^2+y^2\,=\,1$ 
and a 0-dimensional isolated point $(-2,3)$.
Starting from $(x_0,y_0)\,=\,(1.8,0.6)$,
the rank-1 Newton's iteration
(\ref{ni}) converges $(\hat{x},\hat{y})\,\approx\,
(0.928428592, 0.3715109)$ on the unit circle.
%
Starting from $(x_0,y_0)\,=\,(0.4,0.2)$, the iteration converges to
the point
$(0.8007609\ldots,\,0.5989721\ldots)$.
Both sequences of iterates asymptotically follow corresponding normal lines 
of the solution set toward the particular solutions as shown in 
Figure~\ref{f:dir}. 
\end{example}

\section{Algorithms for computing 
\mbox{$J_\rkr(\bdx_k)^\dagger\,\bdf(\bdx_k)$}}

Computing the iteration shift $J_\rkr(\bdx_k)^\dagger\,\bdf(\bdx_k)$ at the 
iterate
$\bdx_k$ is the problem of calculating the {\em minimum norm solution}, or
{\em minimum norm least squares solution} if $\bdb\,\not\in\,\cR(A_\rkr)$, 
of the linear system
\begin{equation}\label{Arzb}  A_\rkr\,\bdz ~~=~~ \bdb
~~~~\mbox{where}~~~~ A\,\in\,\C^{m\times n}
~~~~\mbox{and}~~~~ r\,\le\,\min\{m,\,n\}.
\end{equation}
The most reliable method is based on the singular value
decomposition in the following Algorithm 1.

\vspace{2mm}
\begin{itemize}
\item[] {\bf Algorithm 1: Computing $A_\rkr^\dagger\,\bdb$ by SVD}
\item[] Input: matrix $A\,\in\,\C^{m\times n}$, vector $\bdb\,\in\,\C^m$,
integer $r\,\le\,\min\{m,\,n\}$.
\begin{itemize}
\item By a full or partial SVD, calculate singular values and corresponding
left and right singular vectors 
$\sg_j,\,\bdu_j,\,\bdv_j$, for $j\,=\,1,\ldots,r$.
\item calculate $\bdy\,=\,\left[\frac{\bdu_1^\h\,\bdb}{\sg_1},\ldots,
\frac{\bdu_r^\h\,\bdb}{\sg_r}\right]^\top$ 
\item calculate $\bdz\,=\,[\bdv_1,\ldots,\bdv_r]\,\bdy$
\end{itemize}
\item[] output $A_\rkr^\dagger\,\bdb\,=\,\bdz$
\end{itemize}

\vspace{2mm}
If the Jacobian $J(\bdx_*)$ at the zero $\bdx_*$ is of low rank
such that $r\,\ll\,\min\{m,\,n\}$, a partial SVD such as the
USV-plus decomposition \cite{lee-li-zeng} can be more efficient.
For the cases of $r\,\lessapprox\,\min\{m,n\}$, we need the following lemma
\cite[Lemma 13]{Zeng2019}.

\begin{lemma} \label{l:mns}
Let $A$ be an $m\times n$ matrix with right singular vectors
$\bdv_1,\,\ldots,\bdv_n$ and $A_\rkr$ be its rank-$r$ projection. 
Assume $\sg_r(A)\,>\,\sg_{r+1}(A)$ and $N$ is a matrix whose columns
form an orthonormal basis for $\spn\{\bdv_{r+1},\ldots,\bdv_n\}$.
Then, for every $m$-dimensional vector $\bdb$, the following identity hold:
\begin{equation}\label{Ardb}  A_\rkr^\dagger\,\bdb ~~=~~ (I-N\,N^\h)\,
\left[\begin{array}{c}\mu\,N^\h\\ A \end{array}\right]^\dagger\,
\left[\begin{array}{c}\bdo \\ \bdb \end{array}\right].
\end{equation}
\end{lemma}

For matrices $A\,\in\,\C^{m\times n}$ of large sizes with 
$r\,\approx\,\min\{m,\,n\}$, the SVD can be unnecessarily expensive.
There are reliable algorithms that can be substantially more efficient.
We briefly elaborate some of those algorithms in this section.

The simplest case of computing $A_\rkr^\dagger\,\bdb$ is when $r$ equals 
the row dimension $n$.
This case arises when solving an underdetermined system 
$\bdf(\bdx)\,=\,\bdo$ when the Jacobian
is surjective at the nonisolated solution. 
The following algorithm is a well established numerical
methods for computing the minimum norm solution of a full rank underdetermined
linear system \cite[\S5.6.2]{golub-vanloan4}.

\vspace{2mm}
\begin{itemize}
\item[] {\bf Algorithm 2: Solving (\ref{Arzb}) when $r\,=\,m\,<\,n$}
\item[] Input: matrix $A\,\in\,\C^{m\times n}$ with $m\,<\,n$, 
vector $\bdb\,\in\,\C^m$ (integer $r\,=\,m$).
\begin{itemize}
\item calculate the thin QR decomposition \cite[p. 248]{golub-vanloan4}
 $A^\h\,=\,Q\,R$
\item solve the lower-triangular system $R^\h\,\bdy\,=\,\bdb$ for $\bdy$
\item set $\bdz\,=\,Q\,\bdy$
\end{itemize}
\item[] output $A_\rkr^\dagger\,\bdb\,=\,\bdz$
\end{itemize}

\vspace{2mm}
On the case $r\,<\,m\,<\,n$, the following Algorithm 3 is modified from
the rank-revealing method in \cite{li-zeng-03}.

\vspace{2mm}
\begin{itemize}
\item[] {\bf Algorithm 3: Solving (\ref{Arzb}) when $r\,<\,m\,<\,n$} 
\item[] Input: matrix $A\,\in\,\C^{m\times n}$ with $m\,<\,n$, 
vector $\bdb\,\in\,\C^m$, integer $r\,<\,m$.
\begin{itemize}
\item calculate the thin QR decomposition \cite[p. 248]{golub-vanloan4}
 $A^\h\,=\,Q_0\,R_0$ 
\item set $G_0\,=\,R_0^\h$.
\item for $k\,=\,1,\ldots,m-r$ do
\begin{itemize}
\item calculate the vector $\bdu_k$ as the terminating iterate of the 
iteration (see \cite{li-zeng-03} for details) with $\tau\,=\,\|A\|_\infty$
\begin{equation}\label{lzit}
\bdy_{j+1} = \bdy_j-
\left[\begin{array}{c} 2\tau\,\bdy_j^\h \\ G_{k-1} 
\end{array}\right]^\dagger 
\left[\begin{array}{c} \tau\,\bdy_j^\h\,\bdy_j-\tau \\ G_{k-1}\,\bdy_j 
\end{array}\right], ~~~j=0,1,\ldots
\end{equation}
from a random unit vector $\bdy_0\,\in\,\C^m$
\item As a by product of terminating (\ref{lzit}), extract the thin QR 
decomposition $[2\,\tau\,\bdu_k, \, G_{k-1}^\h]^\h \,=\,Q_k\,G_k$ 
\end{itemize}
\item[] end do
\item solve for $\bdy\,=\,\bdy_*$ of the triangular system 
\[ G_{m-r}^\h\,\bdy\,=\,Q_{m-r}^\h\left[\begin{array}{l} \bdo_{m-r} \\ \bdb 
\end{array}\right]
\]
\item set $\bdz_*\,=\,Q_0\,(I-U\,U^\h)\,\bdy_*$ where $U\,=\,[\bdu_1,\ldots,
\bdu_{m-r}]$
\end{itemize}
\item[] output $A_\rkr^\dagger\,\bdb\,=\,\bdz_*$
\end{itemize}

\vspace{2mm}
Notice that $G_0$ is lower triangular.
It is a standard technique in numerical linear algebra to apply Given's 
rotation \cite[\S5.1.8]{golub-vanloan4} to obtain a QR decomposition
$[2\,\tau\,\bdy_j, \, G_{k-1}^\h]^\h \,=\,Q\,R$ where $R$ is in lower
triangular form throughout the process. 
The algorithm can be easily explained as follows: Let 
$A^\h\,=\,[Q_0,\tilde{Q}]\,\left[\begin{array}{c} R_0 \\ O \end{array}\right]$
be a full QR decomposition of $A^\h$ 
(while $Q_0\,R_0$ is the corresponding thin version).
Then $\cR(\tilde{Q})$ is a subspace of $\cK(A)$ whose basis consists of
$m-r$ additional vectors besides columns of $\tilde{Q}$.
The vectors $Q_0\,\bdu_1,\ldots,Q_0\,\bdu_{m-r}$ approximately form an 
orthonormal basis for the vector space spanned by the $m-r$ right singular 
vectors of $A$. 
The solution $\bdy_*$ of the equation 
$G_{m-r}^\h\,\bdy\,=\,Q_{m-r}^\h\bdb$ is the least squares solution of
the linear system
\[ \left[\begin{array}{c} 2\,\tau\,U^\h \\ R_0^\h \end{array}\right]\,\bdy
 ~~=~~ \left[\begin{array}{c} \bdo \\ \bdb \end{array}\right].
\]
Then it is a straightforward verification using Lemma~\ref{l:mns} that 
$\bdz_*\,=\,Q_0\,(I-U\,U^\h)\,\bdy_*$ is the least squares solution of
$A_\rkr\,\bdz\,=\,\bdb$ that is orthogonal to columns of 
$[Q_0\,U,\,\tilde{Q}]$.

For the cases of solving (\ref{Arzb}) where $m\,\ge\,n$, the kernel or the
partial kernel of $A$ generally can not be computed from a QR decomposition
of $A^\h$ like the cases for $m\,<\,n$.
As a result, a basis $\{\bdu_1,\,\ldots,\,\bdu_{n-r}\}$ for $\cK(A_\rkr)$
needs to be computed as shown in Algrithm 4 below.

\vspace{2mm}
\begin{itemize}
\item[] {\bf Algorithm 4:Solving (\ref{Arzb}) when $r\,\le\,n\,\le\,m$} 
\item[] Input: matrix $A\,\in\,\C^{m\times n}$ with $m\,\ge\,n$, 
vector $\bdb\,\in\,\C^m$, integer $r\,\le\,n$.
\begin{itemize}
\item calculate the thin QR decomposition \cite[p. 248]{golub-vanloan4}
 $A\,=\,Q_0\,R_0$ 
\item for $k\,=\,1,\ldots,n-r$ do
\begin{itemize}\parskip-1mm
\item calculate the vector $\bdu_k$ as the terminating iterate of the 
iteration (see \cite{li-zeng-03} for details) with $\tau\,=\,\|A\|_\infty$
\begin{equation}\label{lzit1}
\bdy_{j+1} = \bdy_j-
\left[\begin{array}{c} 2\tau\,\bdy_j^\h \\ R_{k-1} 
\end{array}\right]^\dagger 
\left[\begin{array}{c} \tau\,\bdy_j^\h\,\bdy_j-\tau \\ R_{k-1}\,\bdy_j 
\end{array}\right], ~~~j=0,1,\ldots
\end{equation}
from a random unit vector $\bdy_0\,\in\,\C^m$
\item As a by product of terminating (\ref{lzit}), extract the thin QR 
decomposition $[2\,\tau\,\bdu_k, \, R_{k-1}]^\h \,=\,Q_k\,R_k$ 
\end{itemize}
\item[] end do
\item solve for $\bdy\,=\,\bdy_*$ of the triangular system 
\[ R_{n-r}\,\bdy\,=\,Q_{n-r}^\h\left[\begin{array}{l} \bdo_{n-r} \\ \bdb 
\end{array}\right]
\]
\item set $\bdz_*\,=\,(I-U\,U^\h)\,\bdy_*$ where $U\,=\,[\bdu_1,\ldots,
\bdu_{n-r}]$
\end{itemize}
\item[] output $A_\rkr^\dagger\,\bdb\,=\,\bdz_*$
\end{itemize}

\section{Modeling with nonisolated solutions}\label{s:m}

\begin{wrapfigure}{r}{3.2in}
\begin{center}
\epsfig{figure=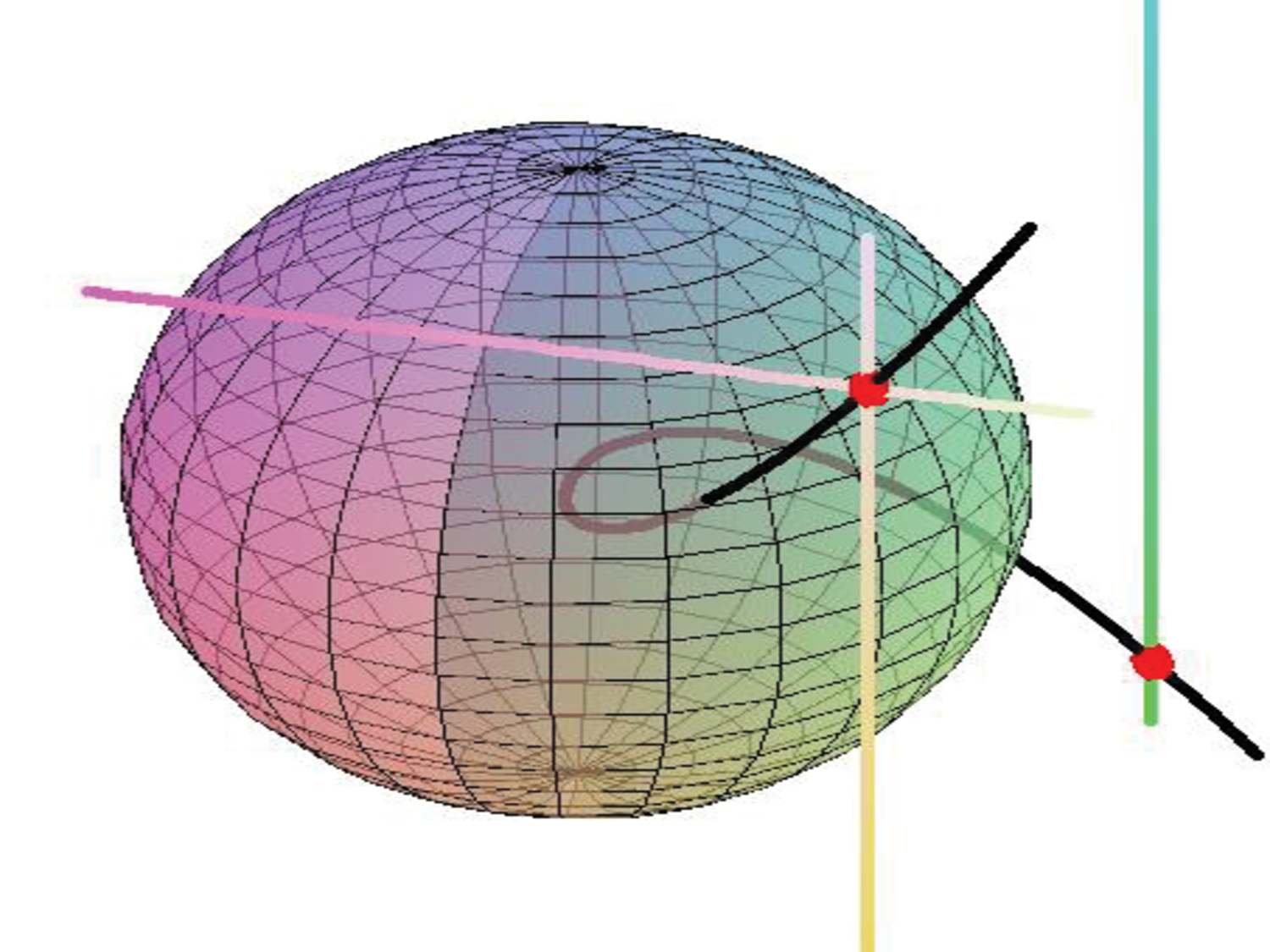,width=3.15in,height=3.0in}
\end{center} \vspace{-8mm}
\caption{\small Solution sets of the system 
in (\ref{ilu})}\label{f:ilex}
\end{wrapfigure}
Models with singular nonisolated solutions arise in many applications.
Unaware of the capability of Newton's iteration in computing such solutions,
scientific computing practitioners go to great lengths to make 
solutions isolated by various techniques such as adding auxiliary equations
and variables. 
We shall elaborate case studies in which nonisolated solutions
can naturally be modeled into well-posed computational problems.
The dimensions of the solution sets and their semiregularity may also be 
obtained analytically in the modeling process so that the projection rank $r$
of Newton's iteration (\ref{ni}) becomes readily available.
The rank-$r$ Newton's iteration for any integer $r$ is implemented in our
Numerical Algebraic Computation toolbox {\sc NAClab}
for Matlab as a functionality {\tt Newton}.
We shall demonstrate the implementation and its effectiveness in 
computing nonisolated solutions with no need for auxiliary equations.

\subsection{Numerical Algebraic Geometry}\label{s:nag}

One of the main subjects of numerical algebraic geometry and its application
in algebraic kinematics is computing solutions of positive dimensions to 
polynomial systems, as elaborated extensively by Wampler and Sommese 
\cite{WamSom11} and in the monographs \cite{BSHW13,som-wam}.
Various mechanisms have been developed in solving those polynomial systems 
for nonisolated solutions, including adding auxiliary equations and variables
to isolate witness points on the solution sets. 
The rank-$r$ Newton's iteration (\ref{ni}) can be applied to calculating 
semiregular witness points directly on the polynomial systems without needing 
extra equations and variables.

\begin{example}[A polynomial system]\label{e:e2}\em
An illustrative example for nonisolated solutions of different dimensions
is given in \cite[p. 143]{BSHW13} as follows:
\begin{equation}\label{ilu} \bdf(x,y,z) = \left(\begin{array}{c}
(y-x^2)\,(x^2+y^2+z^2-1)\,(x-1) \\
(z-x^3)\,(x^2+y^2+z^2-1)\,(y-1) \\
(y-x^2)\,(z-x^3)\,(x^2+y^2+z^2-1)\,(z-1) \end{array}\right).
\end{equation}
Among the solution sets, a point $(1,1,1)$, 
a curve $\{y\,=\,x^2,\,z\,=\,x^3\}$ along with three lines, 
and a surface $\{x^2+y^2+z^2\,=\,1\}$ are of dimensions $0$, $1$ and 
$2$ respectively, as shown in 
Figure~\ref{f:ilex}.
Zeros are semiregular except at intersections of the solution sets.
The iteration (\ref{ni}) converges at quadratic rate toward
those solutions from proper initial iterates by setting the projection 
rank $r\,=\,3$, $2$ and $1$ respectively.
\end{example}

\begin{example}[Perturbed cyclic-4 system]\label{e:cyc}\em
Cyclic-$n$ roots are among the benchmark problems in solving polynomial 
systems. 
Those systems arise in applications such as biunimodular
vectors, a notion traces back to Gauss \cite{FuhRze}.
In general, every cyclic-$n$ system possesses positive dimensional solution
sets if $n$ is a multiple of a perfect square \cite{Backelin}.

We simulate practical computation with empirical data through
the cyclic-4 system in $\bdx\,=\,(x_1,x_2,x_3,x_4)\in\C^4$ with a 
parameter $t\in\C$:
\begin{equation}\label{f4}  \bdf(\bdx,\,t) ~:=~ \left(\begin{array}{c}
x_1+x_2+x_3+x_4 \\
t\,x_1\,x_2+x_2\,x_3+x_3\,x_4+x_4\,x_1 \\
x_1\,x_2\,x_3+x_2\,x_3\,x_4+x_3\,x_4\,x_1+x_4\,x_1\,x_2 \\
x_1\,x_2\,x_3\,x_4-1 \end{array}\right)
\end{equation}
With $t_*\,=\,1$, the equation $\bdf(\bdx,\,1)\,=\,\bdo$ is
the cyclic-4 system whose solution consists of two 1-dimensional sets
\begin{equation}\label{c4z}
 \{x_1\,=\,-x_3,\,x_2\,=\,-x_4,\,x_3\,x_4\,=\,\pm 1,~t\,=\,1\}.
\end{equation}
All zeros are semiregular except eight ultrasingular zeros in the form
of $(\mbox{\footnotesize $\pm 1,\pm 1,\pm 1,\pm 1$})$ 
and $(\mbox{\footnotesize $\pm i,\pm i,\pm i,\pm i$})$ with 
proper choices of signs.
When the parameter $t$ is perturbed from $t_*\,=\,1$ to any other 
nearby value
$\tilde{t}$, the 1-dimensional solution sets dissipate into 16 isolated
solutions. 
Consequently, the parameter value $t_*\,=\,1$ is a bifurcation
point at which the solution changes structure. 

For instance, suppose the parameter $t_*$ is known approximately, say 
$\tilde{t}\,=\,\mbox{\scriptsize 0.9999}$. 
Even though the solutions of $\bdf(\bdx,\tilde{t})\,=\,\bdo$ are all
isolated, Theorem~\ref{t:mt2} ensures the iteration 
\[  \bdx_{k+1} ~~=~~ \bdx_k - 
\bdf_\bdx\big(\bdx_k,\tilde{t}\,\big)_\rk{3}^\dagger\,\,\bdf(\bdx_k,\tilde{t}),
~~~~k\,=\,0,1,\ldots
\]
converges to a numerical solution $\tilde\bdx$ as a zero of the mapping
$\bdx\mapsto\bdf_\bdx(\bdx,\tilde{t})_\rk{3}^\dagger
\bdf(\bdx,\tilde{t})$ that approximates a point in the 1-dimensional
solution set of the underlying system 
$\bdf(\bdx,1)\,=\,\bdo$ with an error in the order of 
$|\tilde{t}-t_*|\,=\,10^{-4}$ if $\bdx_0$ is sufficiently close to
a semiregular point in the solution set (\ref{c4z}), as demonstrated
in the following calling sequence applying the {\sc NAClab} module 
{\tt Newton}:

\setstretch{0.75}{\scriptsize \noindent
\newline $~~~${\verb|>> P = {'x1+x2+x3+x4','0.9999*x1*x2+x2*x3+x3*x4+x4*x1',...;|}~~\blue{\tt \% enter the perturbed cyclic-4 system as}
\newline $~~~${\verb|       'x1*x2*x3+x2*x3*x4+x3*x4*x1+x4*x1*x2', 'x1*x2*x3*x4-1'}; |}~~\blue{\tt \%       a cell array of character strings}
\newline $~~~${\verb|>> v = {'x1','x2','x3','x4'};                              |}~~\blue{\tt \% enter cell array of the variable names}
\newline $~~~${\verb|>> J = PolynomialJacobian(P,v);                     |}~~\blue{\tt \% Jacobian of P w.r.t. the variable names in v}
\newline $~~~${\verb|>> f = @(x,P,J,v) PolynomialEvaluate(P,v,x);   |}~~\blue{\tt \% the function handle for evaluate the system P at x}
\newline $~~~${\verb|>> fjac = @(x,x0,P,J,v) PolynomialEvaluate(J,v,x0)*x;              |}~~\blue{\tt \% function for evaluating J at x}
\newline $~~~${\verb|>> domain = ones(4,1);  parameter = {P,J,v};         |}~~\blue{\tt \% domain (space of 4x1 vectors) and parameters}
\newline $~~~${\verb|>> z0 = [0.8;1.2;-0.8;-1.2];                                                           |}~~\blue{\tt \% initial z0}
\newline $~~~${\verb|>> [z,res,fcond] =  Newton({f,domain,parameter},{fjac,3},z0,1); |}~~\blue{\tt \% call rank-3 Newton iteration on f} 
\newline $~~~~${\verb|                                                         |}~~\blue{\tt \%   projection from z0 using display type 1}
\newline $~~~~${\verb|Step  0:  residual = 7.8e-02  |}
\newline $~~~~${\verb|Step  1:  residual = 2.4e-03  shift = 2.4e-02 |}
\newline $~~~~${\verb|Step  2:  residual = 1.0e-04  shift = 6.8e-04 |}
\newline $~~~~${\verb|Step  3:  residual = 1.0e-04  shift = 5.8e-07 |}
\newline $~~~~${\verb|Step  4:  residual = 1.0e-04  shift = 4.3e-13 |}
\newline $~~~~${\verb|Step  5:  residual = 1.0e-04  shift = 3.6e-16 |}
\newline $~~~~${\verb|Step  6:  residual = 1.0e-04  shift = 1.5e-16 |}
}

\vspace{2mm}
\setstretch{1.0}{}
The iteration does {\em not} converges to a solution to 
$\bdf(\bdx,0.9999)\,=\,\bdo$ as the residual 
$\|\bdf(\bdx_j,0.9999)\|_2$ can only be reduced to $10^{-4}$
but the shifts 
\[ \|\bdx_{j+1}-\bdx_j\|_2 ~=~
\|\bdf_\bdx(\bdx_j,0.9999)_\rk{3}^\dagger\,\bdf(\bdx_j,0.9999)\|_2
~\longrightarrow~ 1.51\times 10^{-16}
\]
approaches to the unit roundoff, indicating the iteration converges to
a stationary point $\tilde\bdx$ at which 
$\bdf_\bdx(\tilde\bdx,0.9999)_\rk{3}^\dagger\,\bdf(\tilde\bdx,0.9999)\,=\,
\bdo$. 
The module {\tt Newton} terminates at $\tilde\bdx$ that 
approximates the nearest solution $\hat\bdx$ of the underlying 
equation $\bdf(\bdx,1)\,=\,\bdo$ where
\begin{align*}
\tilde\bdx &~=~ (\mbox{\scriptsize \tt 0.82287906\,1867739, 
1.21524540\,1950727, -0.82287906\,2858240, -1.215245403\,413521}) \\
\hat\bdx &~=~ (\mbox{\scriptsize \tt 0.82287906\,3773473, 1.21524540\,3637205, 
-0.82287906\,3773473, -1.215245403\,637205}) 
\end{align*}
with a forward error $2.71\times 10^{-9}$ much smaller than the data 
error $10^{-4} = |\tilde{t}-t_*|$ as asserted in Theorem~\ref{t:mt2}.
As a confirmation of the geometric interpretation of the iteration elaborated
in \S\ref{s:geo}, the point $\hat\bdx$ is also only about 
$6.3\times 10^{-4}$ away from the nearest point in the
solution set to the initial iterate
$\bdx_0\,=\,(\mbox{\scriptsize \tt 0.8, 1.2, -0.9, -1.2})$.
\end{example}

\begin{example}[A bifurcation model]\label{e:bif}\em
\,\,\,\,The parameterized cyclic-4 system (\ref{f4}) leads to a bifurcation model
in the form of the equation
\[  \bdf(\bdx,\,t) ~~=~~ \bdo ~~~~\mbox{for}~~~ (\bdx,t)\,\in\,\R^4\times\R
\]
depending on the parameter $t$.
The solution of $\bdf(\bdx,t)\,=\,\bdo$ for 
$(\bdx,t)\in\,\R^4\times\R$ consists of 16 branches 
$\bdx\,=\,\psi_j(t)$ for $j\,=\,1,\ldots,16$ and $t\in\R$. 
On the plane $t\,=\,1$, eight pairs of those 16 branches intersect at 
8 bifurcation points that are embeded in two additional branches of solution 
curves (\ref{c4z}).
All solutions are semiregular except at the intersections. 
The parameter value $t_*\,=\,1$ is the bifurcation point at which
the solution changes the structure.

After finding $\tilde\bdx$ in Example~\ref{e:e2} as an approximate
zero of $\bdx\,\mapsto\,\bdf(\bdx,\,t_*)$
from the empirical data $\tilde{t}\,=\,0.9999$, the point 
$(\tilde\bdx,\tilde{t})$ is very close to a branch in (\ref{c4z}).
The point $(\tilde\bdx,\tilde{t})$ is likely to be closer to a solution
branch in (\ref{c4z}) than to other 17 branches.
If so, the iteration
\[  (\bdx_{j+1},\,t_{j+1}) ~~=~~ (\bdx_j,\,t_j) - 
\bdf_{\bdx t}(\bdx_j,t_j)_\rk{4}^\dagger~\bdf(\bdx_j,t_j), ~~~~j\,=\,0,1,\ldots
\]
starting from the initial iterate 
$(\bdx_0,\,t_0)\,=\,(\tilde\bdx,\tilde{t})$ converges roughly to the nearest
point in the branch in (\ref{c4z}) on the plane $t\,=\,1$ in $\R^4\times\R$
by the geometric interpretation elaborated in \S\ref{s:geo}. 
This is significant because the sequence $\{t_j\}$ converges to the 
important bifurcation point $t_*\,=\,1$.
This observation can easily be confirmed by {\sc NAClab} {\tt Newton} and
obtains a solution $(x_1,x_2,x_3,x_4,\,t)$ as
\[
(\mbox{\scriptsize \tt 
0.822879063773473, 1.215245403637205, -0.822879063773474, -1.215245403637204,
1.000000000000000}) 
\]
that is accurate with at least 15 digits, particularly the bifurcation point
at $t\,=\,1.0$.
The residual $\big\|\bdf(\bdx_j,t_j)\big\|_2$ reduces from $10^{-4}$ 
to $4.44\times 10^{-16}$, indicating convergence to a solution of
$\bdf(\bdx,t)\,=\,\bdo$.
Also notice the $\bdx$ component is identical to nearest point
$\hat\bdx$ for at least 15 digits, confirming the geometric interpretation 
in \S\ref{s:geo} again.
\end{example}

\begin{remark}[Comparison with existing methods]\label{r:slc}\em
Regularization has been required for computing nonisolated solutions.
An established strategy for regularizing solutions of positive dimensions is 
called {\em linear slicing}~\cite[\S 13.2]{som-wam},
namely augmenting the system with auxiliary linear equations. 
Every auxiliary equation generically reduces the dimension by one. 
The resulting overdetermined system $\bdf(\bdx)\,=\,\bdo$ can then 
be squared by {\em randomization} in the form of 
$\Lambda\,\bdf(\bdx)\,=\,\bdo$ where 
$\Lambda$ is a random matrix of proper size \cite[\S 13.5]{som-wam} so that
conventional Newton's iteration (\ref{ni0}) applies.
Besides creating extraneous solutions by randomization, two issues arise:

\begin{itemize}
\item The singular solutions may disappear when the system is given with 
empirical data as we simulated in Example~\ref{e:cyc}. 
The augmented slicing equations may not intersect a solution. 
\item Does the randomized square system with perturbed data has zeros 
approximate the solutions of the underlying system with exact data?
\end{itemize}

Assuming the linear slicing isolates regular solutions, Theorem~\ref{t:mt2} 
actually implies that direct applying the Gauss-Newton iteration on the 
regularized system without randomization results in local linear convergence 
to a stationary point approximating an underlying solution even if the 
solution set is annihilated by the data perturbation. 
The randomization likely achieves a similar result.
From our experiment, linear slicing combined with randomization produces end
results of near identical accuracy with a moderate increase of condition 
numbers.
The advantage of the rank-$r$ Newton's iteration is to eliminate the need of
the solution isolation and randomization processes and to enable direct 
convergence to the solution manifolds.
\end{remark}

\subsection{Numerical greatest common divisor}\label{s:gcd}

For a polynomial pair $p$ and $q$ of degrees $m$ and
$n$ respectively with a greatest common divisor (GCD) $u$ of degree $k$
along with cofactors $v$ and $w$, the GCD problem can be modeled as a 
zero-finding problem
\begin{align}\label{gcdeq}
\bdf(u,v,w) &~~=~~ (0,0) ~\in~ \mP_m\times\mP_n \\
&~~~~\mbox{for}~~~ (u,v,w)~\in~ \mP_k\times\mP_{m-k}\times\mP_{n-k} 
\nonumber
\end{align}
with the holomorphic mapping
\begin{equation}\label{fgcd}
  \mapform{\bdf}{\mP_k\times\mP_{m-k}\times\mP_{n-k}}{\mP_m\times\mP_n}{
(u,v,w)}{(u\,v-p,~u\,w-q)}
\end{equation}
where $\mP_m$ is the vector space of polynomials with degrees up to $m$, 
etc.
Through standard bases of monomials, the mapping $\bdf$ in (\ref{fgcd}) 
can be represented by a mapping from $\C^{m+n-k+3}$ to $\C^{m+n+2}$.
Let $(u_*,\,v_*,\,w_*)$ denote a particular solution of (\ref{gcdeq}).
The general solution of (\ref{gcdeq}) is a 1-dimensional set
\begin{equation}\label{gs}
  \left\{\left. \big(t\,u_*,\,\mbox{$\frac{1}{t}$}\,v_*,\,
\mbox{$\frac{1}{t}$}\,w_*\big) ~\right|~ t\,\in\,\C\setminus\{0\}\right\}
\end{equation}
on which the Jacobian is rank-deficient. 
Adding one auxiliary equation, however, the Jacobian of the modified
mapping becomes injective \cite{ZengGCD}, implying the Jacobian of $\bdf$ 
in (\ref{fgcd}) is of rank deficient by one columnwise.
As a result, the set (\ref{gs}) consists of semiregular zeros of the mapping 
$\bdf$ in (\ref{fgcd}) and the iteration (\ref{ni}) is locally quadratically 
convergent by setting the projection rank
\[ r~=~(k+1)+(m-k+1)+(n-k+1)-1 ~=~ m+n-k+2.
\]
Extra equations are not needed.
Although there are infinitely many solutions forming a 1-dimensional set, 
there is practically no difference in finding anyone over the other. 

\begin{example}[Numerical greatest common divisor]\label{e:gcd}\em
For a simple demo of the iteration (\ref{ni}) and its
{\sc NAClab} implementation {\tt Newton}, let
\begin{equation}\label{pq}
 p ~=~ 
\mbox{\scriptsize $-1.3333$} \mbox{\scriptsize $ -2.3333$}\,x 
\mbox{\scriptsize $-4$}\,x^2 \mbox{\scriptsize $-3.6667$}\,x^3 
\mbox{\scriptsize $-2.6667$}\,x^4 - x^5,
~~~\mbox{and}~~~ q ~=~ 
\mbox{\scriptsize $-1.9999$} + x + x^2 \mbox{\scriptsize $+ 3$}\,x^3
\end{equation}
that are considered perturbed data of an underlying polynomial pair with a GCD 
$1+x+x^2$.
The process of identifying the GCD degree and computing the initial
approximations of the GCD along with the cofactors can be found in
\cite{ZengGCD}. 
With {\sc NAClab} installed, the following sequence of Matlab statements 
carries out the rank-8 Newton's iteration
\[  (u_{j+1},\,v_{j+1},\,w_{j+1}) ~~=~~ (u_j,\,v_j,\,w_j) -
J_\rk{8}(u_j,\,v_j,\,w_j)^\dagger\,\bdf(u_j,\,v_j,\,w_j)
\]
for $j\,=\,0,1,\ldots$.
\end{example}

\setstretch{0.75}{\scriptsize \noindent
$~~${\verb|>> p = '-1.3333-2.3333*x-4*x^2-3.6667*x^3-2.6667*x^4-x^5';         |}~~\blue{\tt \% enter p as a character string}
\newline $~~~${\verb|>> q = '-1.9999+x+x^2+3*x^3';                                                  |}~~\blue{\tt \% enter q similarly}
\newline $~~~${\verb|>> f = ...                    |}~~\blue{\tt \% enter the mapping as function handle for (u,v,w) -> (u*v-p, u*w-q)}
\newline $~~~${\verb|       @(u,v,w,p,q) {pminus(ptimes(u,v),p),pminus(ptimes(u,w),q)};      |}
\newline $~~~${\verb|>> J = ...   |}~~\blue{\tt \% enter the Jacobian J at (u0,v0,w0) as the mapping (u,v,w) -> (u0*v+u*v0, u0*w+u*w0)}
\newline $~~~${\verb|        @(u,v,w,u0,v0,w0,p,q){pplus(ptimes(u0,v),ptimes(u,v0)),pplus(ptimes(u0,w),ptimes(u,w0))}; |}
\newline $~~~${\verb|>> domain = {'1+x+x^2','1+x+x^2+x^3', '1+x'};     |}~~\blue{\tt \% representation of the domain for the mapping f}
\newline $~~~${\verb|>> parameter = {p,q};                                               |}~~\blue{\tt \% parameters for the mapping f}
\newline $~~~${\verb|>> u0 = 'x^2+1.4*x+1.6';  v0 = '-1.5-x-1.6*x^2-x^3'; w0 = '-2+2.8*x';         |}~~\blue{\tt \% initial (u0,v0,w0)}
\newline $~~~${\verb|>> [z,res,fcond] = Newton({f,domain,parameter},{J,8},{u0,v0,w0},1);        |}~~\blue{\tt \% Newton on f with J in}
\newline $~~~${\verb|                    |}~~\blue{\tt \%  rank-8 projection from the initial iterate (u0,v0,w0) with display setting 2}
\newline $~${\verb|Step  0:  residual = 1.5e+00                   1.6+14.*x+x^2|}
\newline $~${\verb|Step  1:  residual = 1.1e-01  shift = 4.9e-01  1.114771189108 + 1.114523702576*x + 1.088690420621*x^2|}
\newline $~${\verb|Step  2:  residual = 1.2e-03  shift = 5.9e-02  1.089839864820 + 1.090023690710*x + 1.089788605779*x^2 |}
\newline $~${\verb|Step  3:  residual = 8.4e-06  shift = 1.0e-03  1.089756319215 + 1.089767203330*x + 1.089783432095*x^2 |}
\newline $~${\verb|Step  4:  residual = 8.3e-06  shift = 1.4e-07  1.089756333892 + 1.089767171466*x + 1.089783428226*x^2 |}
\newline $~${\verb|Step  5:  residual = 8.3e-06  shift = 5.1e-13  1.089756333892 + 1.089767171466*x + 1.089783428226*x^2|}
\newline $~${\verb|Step  6:  residual = 8.3e-06  shift = 7.0e-15  1.089756333892 + 1.089767171469*x + 1.089783428226*x^2|}
}

\vspace{2mm}
\setstretch{1.0}{ }
The computed GCD 
$\mbox{\scriptsize $1.08976$} + \mbox{\scriptsize $1.08977$}\,x + 
\mbox{\scriptsize $1.08978$}\,x^2$ is of 
the scaling independent distance $1.02\times 10^{-5}$ that is in the
same order of the data error $2.41\times 10^{-5}$.

Similar to Example~\ref{e:bif}, the system with inexact data for $(p,q)$ 
does not have a solution as
the residual $\|\bdf(u_j,v_j,w_j)\|_2$ can not be reduced below $8.0\times
10^{-6}$ but the shift 
$\|(u_{j+1},v_{j+1},w_{j+1})-(u_j,v_j,w_j)\|_2$ reduces to near 
unit roundoff, implying
$J_\rk{8}(u_j,v_j,w_j)^\dagger\,\bdf(u_j,v_j,w_j)$ approaches zero.

\begin{remark}[On existing methods]\em
The 1-dimensional solution (\ref{gs}) can be isolated by linear slicing
(c.f. Remark~\ref{r:slc}) such as requiring the GCD to be monic or 
to satisfy a linear equation \cite{ZengGCD}.
The rank-$r$ Newton's iteration achieves the results of the same or better
accuracy
in practical computation without needing an isolation process.
\end{remark}

\subsection{Accurate computation of defective eigenvalues}\label{s:eig}

Accurate computation of defective eigenvalues requires regularization due to 
hypersensitivities to data perturbations \cite{pseudoeig}.
Let $\hat\la\,\in\,\C$ be an eigenvalue of $A\,\in\,\C^{n\times n}$ with 
what we call a multiplicity support $m\times k$. 
Namely $\hat\la$ is of the geometric multiplicity $m$ with the smallest 
Jordan block of size $k\times k$.
The problem of computing $\hat\la$ can be naturally modeled as solving 
the equation
\begin{align}\label{eigeq}
\bdg(\la,\,X,\,A) &~~=~~ O ~\in~ \C^{n\times k} \\
&~~~~~\mbox{for}~~~ (\la,X) ~\in~ \C\times\C^{n\times k} \nonumber
\end{align}
where
\begin{equation}\label{geig}
\mapform{\bdg}{\C\times\C^{n\times k}\times\C^{n\times n}}{\C^{n\times k}}{
(\la,\,X,\,G)}{(G-\la\,I)\,X-X\,S} 
\end{equation}
with a constant nilpotent upper triangular matrix parameter 
$S\,\in\,\C^{k\times k}$ of rank $k-1$.
It can be verified with a standard linear algebra that the solution set of
(\ref{eigeq}) is of dimension $\ell\,=\,m\,k$ in the form of 
$(\la_*,\,\al_1\,X_1+\cdots+\al_{\ell}\,X_{\ell})$ with a constant 
eigenvalue component $\la_*$.
By \cite[Lemma~2(ii)]{pseudoeig}, the Jacobian 
$\bdg_{\la X}(\hat\la,\hat{X},A)$ is of nullity $m\,k$,
implying the solution is semiregular.
Similar to the case study in \S\ref{s:gcd}, only a representative in the
solution set is needed.

Theorem~\ref{t:mt} and Theorem~\ref{t:mt2} ensure a defective eigenvalue 
can be accurately computed by the rank-$(n\,k-m\,k+1)$ Newton's
iteration~(\ref{ni}).

\begin{example}[Defective eigenvalue computation]\label{e:eig}\em
We show the effectiveness of the rank-$r$ Newton's iteration in 
modeling and computation of defective eigenvalues with a simple example.
Let the matrix $A$ and a perturbed matrix $E$ be as follows.
\[  
A = \mbox{\tiny $\left[\begin{array}{rrrrrr}
   -1  &   0  &   3 &    0  &   2  &   1 \\
    1  &   1  &  -1 &    1  &   0  &   0 \\
   -2  &  -1  &   4 &    1  &   1  &   0 \\
    3  &  -3  &  -3 &    5  &  -1  &  -1 \\
   -3  &   1  &   3 &   -1  &   5  &   2 \\
    1  &   0  &  -1 &    0  &  -1  &   2
\end{array}\right]$} \;\;\mbox{and}\;\;
E = 10^{-6} \mbox{\tiny $\left[\begin{array}{rrrrrr}
   .1 & -.7 & -.4 & -1.0 &  .2 &  .6 \\
  -.2 &  .1 & -.1 & -.5 &  .5 &  .0 \\
   .3 & -.8 & -.6 & -.1 &  .4 &  .1 \\
  -.5 &  .0 &  .1 &  .7 & -.2 &  .5 \\
  -.2 & -.2 & -.8 & -.7 & -.4 & -.5 \\
  -.2 & -.1 &  .8 & -.5 & -.7 & -.6
\end{array}\right]$}
\]
Matrix $A\,\in\,\C^{6\times 6}$ possesses an exact eigenvalue 
$\la_*\,=\,3$ with the multiplicity support $2\times 2$ 
(c.f. \cite{pseudoeig} for identifying multiplicity supports). 
At the fixed $\la_0\,=\,2.9$, we can apply the {\sc NAClab} module
{\tt LinearSolve} to solve the homogeneous linear system
$\bdg(\la_0,\,X,A)\,=\,O$ with the rank-8 projection for a solution 
$X\,=\,X_0$ and obtain the initial iterate $(\la_0,\,X_0)$. 
With $A$ and $E$ above entered in Matlab, we apply the rank-$r$ Newton's
iteration (\ref{ni}) with $r\,=\,(n-m)\,k+1\,=\,9$ by executing the NAClab 
module {\tt Newton} in the following calling sequence.
\end{example}

\setstretch{0.75}{\scriptsize \noindent
$~~${\verb|>> g = @(lambda,X,G,S) G*X-lambda*X-X*S;             |}~~\blue{\tt \% enter the mapping g as a function handle}
\newline $~~~${\verb|>> J = @(lambda,X,lambda0,X0,G,S) G*X-lambda*X0-lambda0*X-X*S;             |}~~\blue{\tt \% enter the Jacobian}
\newline $~~~${\verb|>> domain = {1,ones(n,k)};                     |}~~\blue{\tt \% representation of the domain for the mapping g}
\newline $~~~${\verb|>> parameter = {A,S};                                             |}~~\blue{\tt \% parameters of the mapping g}
\newline $~~~${\verb|>> [z,res,fcond] = ...          |}~~\blue{\tt \% Call Newton on the mapping g with J in rank-9 projection from} 
\newline $~~~${\verb|   Newton({g,domain,parameter},{J,9},{lambda0,X0},2);       |}~~\blue{\tt \%    (lambda0,X0) using display type 2}
\newline $~~~~${\verb| |}
\newline $~~~~${\verb|Step  0:  residual = 2.5e-02                     2.900000000000000|}
\newline $~~~~${\verb|Step  1:  residual = 4.1e-03  shift = 1.0e-01    3.000493218848026|}
\newline $~~~~${\verb|Step  2:  residual = 3.0e-06  shift = 6.5e-03    2.999999313752690 |}
\newline $~~~~${\verb|Step  3:  residual = 2.6e-12  shift = 4.1e-06    3.000000000001909|}
\newline $~~~~${\verb|Step  4:  residual = 4.4e-16  shift = 3.2e-12    3.000000000000000|}
\newline $~~~~${\verb|Step  5:  residual = 2.2e-16  shift = 5.1e-16    3.000000000000000|}
}

\vspace{2mm}
\setstretch{1.0}{ }
The $\la$ components of the iteration accurately converges to the 
defective eigenvalue $\la_*\,=\,3$ at roughly the quadratic rate with an
accuracy at the order of the unit roundoff.

We simulate practical computation with imperfect empirical data by using the
perturbed matrix $\tilde{A}\,=\,A+E$ in the iteration
\[ \big(\tilde\la_{j+1},\,\tilde{X}_{j+1}\big) ~~=~~ 
\big(\tilde\la_j,\,\tilde{X}_j\big) -
\bdg_{_{\la X}}\big(\tilde\la_j,\,\tilde{X}_j,\,A+E\big)_\rk{9}^\dagger\, 
\bdg\big(\tilde\la_j,\,\tilde{X}_j,\,A+E\big).
\]
The iteration reaches the numerical solution $(\tilde\la,\,\tilde{X})$
where $\tilde\la\,=\,3.00000102$, with a forward accuracy 
$1.02\times 10^{-6}$ about the same as the data perturbation
$\|E\|_2\,\approx\,1.9\times 10^{-6}$.

\begin{example}[Nearest matrix with the defective eigenvalue]\label{e:eig2}\em
\,\,
Furthermore, the iteration (\ref{ni}) can again be applied to refine the 
eigenvalue computation by solving the equation
\begin{equation}\label{eigeq3}
   \bdg(\la,\,X,\,G) ~~=~~ O ~~~~\mbox{for}~~~~ (\la,\,X,\,G)\,\in\,
\C\times\C^{n\times k}\times\C^{n\times n}
\end{equation}
in the specific iteration
\begin{equation}\label{eigit3}
 \big(\hat\la_{j+1},\,\hat{X}_{j+1},\,\hat{A}_{j+1}\big) ~~=~~ 
\big(\hat\la_j,\,\hat{X}_j,\,\hat{A}_j\big) -
\bdg_{_{\la X G}}\big(\hat\la_j,\,\hat{X}_j,\,\hat{A}_j\big)_\rk{12}^\dagger\, 
\bdg\big(\hat\la_j,\,\hat{X}_j,\,\hat{A}_j\big)
\end{equation}
for $j\,=\,0,1,\ldots,$ starting from 
$(\tilde\la_0,\,\tilde{X}_0,\,\hat{A}_0\big)\,=\,
(\tilde\la,\,\tilde{X},\,A+E)$.
The projection rank $r\,=\,n\,k\,=\,12$ because,
by \cite[Lemma 2(ii)]{pseudoeig}, the (full) Jacobian of $\bdg$ is 
surjective, implying the solution is semiregular with dimension 
$(1+n\,k+n^2)-n\,k\,=\,n^2+1$.

In this example, the refinement (\ref{eigit3}) needs only 1 step to reduce 
the residual of the equation (\ref{eigeq3}) from
$2.99\times 10^{-7}$ to $1.17\times 10^{-15}$.
The iteration terminates at 
$\big(\hat\la,\,\hat{X},\,\hat{A}\big)$ with 
$\hat\la\,\approx\,3.000000103$, and the third component $\hat{A}$ is
the matrix with an exact defective eigenvalue $\hat\la$ of multiplicity
support $m\times k$, implying the backward error is 
$\|\tilde{A}-\hat{A}\|\,\approx\,7.59\times 10^{-7}$.
From the geometric interpretation elaborated in \S\ref{s:geo}, the iteration
(\ref{eigit3}) converges to $(\hat\la,\,\hat{X},\,\hat{A})$ that is roughly
the nearest point in the solution set of (\ref{eigeq3}) from the initial
iterate $(\tilde\la,\,\tilde{X},\,A+E)$.
\end{example}

\begin{remark}[On the existing method]\em
The method proposed in \cite{pseudoeig} appears to be the only published method
for accurate computation of defective eigenvalues. 
The method employs \,$m\,k$ \,auxiliary equations
in the eigenvalue model to isolate a 
solution point in the solution set of (\ref{eigeq}).
~Although the computational results are practically equivalent for
Example~\ref{e:eig}, ~those arbitrary extra equations are unnatural and 
are now unnecessary.
Linear slicing for finding the nearest matrix with defective eigenvalues
in Example~\ref{e:eig2} requires quadruple the number of equations
(37 extra equations on top of the 12 existing ones) unnecessarily.
The rank-$r$ Newton's iteration substantially reduces the computing cost
in problems such as Example~\ref{e:eig2}.
\end{remark}

\section{A note on computing ultrasingular zeros}\label{s:sz}

As formulated in \S\ref{s:r}, ultrasingularity occurs if the nullity of the 
Jacobian is higher at a zero than the dimension of the zero. 
Difficulties arise in computing ultrasingular zeros including slow convergence
rate of iterative methods (c.f. \cite{DecKelKel}) and, more importantly, 
barriers of low attainable accuracy \cite{victorpan97,ypma}.
As pointed out in \cite{Kel81}, ``{\em direct} attempts at such computation
may easily fail or give inconclusive results''.
Quadratic convergence rate of Newton's iteration at semiregular 
zeros is not expected at ultrasingular zeros.
At an isolated ultrasingular zero $\bdx_*$ of a smooth mapping 
$\bdf\,:\,\Omega\subset\F^m\rightarrow\F^n$ where $\F\,=\,\R$ or $\C$,
the Jacobian $J(\bdx_*)$ is of rank $r\,<\,m$.
At $k$-th step of conventional Newton's iteration (\ref{ni0}) or the 
Gauss-Newton iteration when $n\,>\,m$, the Jacobian
$J(\bdx_k)$ is usually highly ill-conditioned so that the computation 
of the iterate $\bdx_{k+1}$ is generally inaccurate, substantially limiting
the attainable accuracy of the computed zero even if the iteration converges.

A {\em depth-deflation} strategy \cite{DLZ} deflates the isolated 
ultrasingularity and transforms the zero $\bdx_*$ into a component of a 
possibly regular zero 
$(\bdx_*,\bdy_*)$ of an expanded mapping
\begin{equation}\label{dfm}
\bdg ~~:~~ (\bdx,\bdy)~\mapsto~ \left(\begin{array}{c}
\bdf(\bdx), \, J(\bdx)\,\bdy, \, R\,\bdy - \bde
\end{array}\right)
\end{equation}
where $R$ is a random $(m-r)\times m$ matrix and $\bde\,\ne\,\bdo$.
If $(\bdx_*,\bdy_*)$ is still an ultrasingular zero of 
$(\bdx,\bdy)\,\mapsto\, \bdg(\bdx,\bdy)$, the deflation process can be 
continued recursively.
It is proved in \cite{DLZ} that the number of deflation steps is bounded by
the so-called {\em depth} of an ultrasingular isolated zero. 
In practice, however, one deflation step is likely to be
enough except for the cases where the {\em breadth} 
$\nullity{J(\bdx_*)} = 1$. 
At the terminating step of depth deflation, the ultrasingular zero $\bdx_*$ 
of $\bdf$ is a component of the semiregular zero of the final expanded 
mapping.
As a result, the Gauss-Newton iteration locally converges at quadratic rate.
More importantly, the zero $\bdx_*$ can be computed with an accuracy 
proportional to the data precision or unit round-off, circumventing the 
barrier of the perceived attainable accuracy at ultrasingular zeros.
An earlier deflation strategy in \cite{lvz06} is proven to terminate with
the number of steps bounded by the multiplicity.
A so-called strong deflation method in symbolic-numerical computation is
proposed in \cite{HauWam13} and also proved to terminate in finitely 
many steps.

The deflation strategy applies to systems at ultrasingular nonisolated 
solutions as well. 
We illustrate the deflation process in the following examples.

\begin{example}[Isolated ultrasingularity in a nonisolated zero set]\em
\label{e:is}
\,\,\,\,\,\,\,\, In the
cyclic-4 system in Example~\ref{e:cyc}, the mapping $\bdx\,\mapsto\,
\bdf(\bdx,1)$ as in (\ref{f4}) possesses 8 ultrasingular zeros embedded in the 
two solution curves (\ref{c4z}). 
Those 8 points are nonisolated zeros with isolated ultrasingularity since 
each point is a unique ultrasingular zero in a small open neighborhood.
For instance, the point $\bdx_* = (1,-1,-1,1)$ is such an ultrasingular 
zero at which $\rank{\bdf_\bdx(\bdx_*,1)} = 2$.
Interestingly, for almost all $R\in\R^{2\times 4}$, there is a unique 
$\bdy_*\in\R^4$ such that $(\bdx_*,\bdy_*)$ is a {\em semiregular} isolated
zero of the expanded system $(\bdf(\bdx,1),\bdf_\bdx(\bdx,1)\,\bdy,
R\,\bdy-(1,0))$.
In other words, the deflation method applies to isolated ultrasingularity at 
nonisolated zeros as well and the Gauss-Newton iteration converge 
quadratically to $(\bdx_*,\bdy_*)$.
Even if the system is given with imperfect empirical data, the Gauss-Newton
iteration on the expanded system still converges linearly to a stationary 
point that approximates $(\bdx_*,\bdy_*)$ with an accuracy at the same
order of the data.
The depth deflation method being applicable to nonisolated solutions of
isolated ultrasingularities has apparently not been observed before.
The theoretical termination and the bound on the number of deflation steps 
are still unknown.
\end{example}

\begin{example}[An ultrasingular branch of zeros]\label{e:sb}\em
There are nonisolated zeros \linebreak
where the entire branch is ultrasingular.
Let $\bdx\,=\,(x_1,x_2,x_3,x_4)$ and the mapping
\[ \bdf(\bdx) ~~=~~ \left(\begin{array}{c}
x_1^3+x_2^2+x_3^2\,x_4^2-1 \\ x_1^2+x_2^3+x_3^2\,x_4^2-1 \\
x_1^2+x_2^2+x_3^3\,x_4^3-1
\end{array}\right)
\]
with a 1-dimensional solution curve 
$S\,=\,\{\bdx\,=\,(0,0,s,1/s) ~|~ s\ne 0\}$
on which the Jacobian $J(\bdx)$ is of rank 1.
All the solutions in $S$ are ultrasingular due to 
$\nullity{J(\bdx)} = 3 \ne 1 = \dm_\bdf(\bdx)$ on $S$.
As a result, the equation $\bdf(\bdx) = \bdo$ is underdetermined since,
for every $\bdx\in S$ and almost all matrix $R\in\R^{3\times 4}$, there is
a unique $\bdy\in\R^4$ such that $(\bdx,\bdy)$ is a 1-dimensional zero of
the expanded mapping $(\bdx,\bdy)\,\mapsto\,\bdg(\bdx,\bdy)$ in 
(\ref{dfm}) where $\bde$ can be any nonzero vector, say $(1,0,0)$.
The actual rank-7 Newton's iteration
\[ (\bdx_{k+1},\bdy_{k+1}) ~~=~~ (\bdx_k,\bdy_k) - 
\bdg_{\bdx\bdy}(\bdx_k,\bdy_k)_\rk{7}^\dagger\, \bdg(\bdx_k,\bdy_k),
~~~k\,=\,0,1,\ldots
\] 
can be carried out using the following 
{\sc NAClab} calling sequence.

\vspace{2mm}
\setstretch{0.75}{\scriptsize \noindent
$~~${\verb|>> P = {'x1^3+x2^2+x3^2*x4^2-1'; 'x1^2+x2^3+x3^2*x4^2-1'; 'x1^2+x2^2+x3^3*x4^3-1'};  |}~~\blue{\tt \% enter system}
\newline $~~~${\verb|>> x = {'x1';'x2';'x3';'x4'}; J = pjac(P,x);       |}~~\blue{\tt \% variable name array x, Jacobian of P w.r.t. x}
\newline $~~~${\verb|>> y = ['y1';'y2';'y3';'y4']; R = Srand(3,4);        |}~~\blue{\tt \% expanded variable array y, random 3x4 matrix}
\newline $~~~${\verb|>> F = [P; ptimes(J,y); pminus(ptimes(R,y),[1;0;0])]; K = pjac(F,[x;y]);  |}~~\blue{\tt \% new system and Jacobian}
\newline $~~~${\verb|>> g = @(u,v,F,K,x,y) peval(F,[x;y],[u;v]);           |}~~\blue{\tt \% function handle for expanded mapping g(x,y)}
\newline $~~~${\verb|>> gjac = @(u,v,u0,v0,F,K,x,y) peval(K,[x;y],[u0;v0])*[u;v]; |}~~\blue{\tt \% function handle for expanded Jacobian}
\newline $~~~${\verb|>> u0 = [0.001; 0.003; 0.499; 2.002];                         |}~~\blue{\tt \% an initial estimate of the solution}
\newline $~~~${\verb|>> [~,~,V] = svd(peval(J,v,u0)); v0 = V(:,2:4)*((R*V(:,2:4))\[1;0;0]);         |}~~\blue{\tt \% new initial iterate}
\newline $~~~${\verb|>> [Z,res,fcd] = Newton({g,{ones(4,1),ones(4,1)},{F,K,x,y}}, {gjac,7},{u0,v0}, 1)    |}~~\blue{\tt \% rank-7 Newton}
\newline $~~~~${\verb| |}
\newline $~~~~${\verb|Step  0:  residual = 6.2e-03 |}
\newline $~~~~${\verb|Step  1:  residual = 1.5e-05  shift = 3.0e-03|}
\newline $~~~~${\verb|Step  2:  residual = 1.2e-09  shift = 2.4e-05|}
\newline $~~~~${\verb|Step  3:  residual = 2.2e-16  shift = 1.6e-09|}
\newline $~~~~${\verb|Step  4:  residual = 4.4e-16  shift = 8.2e-16|}
}

\vspace{2mm}
\setstretch{1.0}{ }
The step-by-step {\tt shift} and {\tt residual} show both 
$\|\bdg_{\bdx\bdy}(\bdx_k,\bdy_k)_\rk{7}^\dagger\, \bdg(\bdx_k,\bdy_k)\|_2$
and 
$\|\bdg(\bdx_k,\bdy_k)\|_2$ 
approach zero, implying the limit $(\hat\bdx,\hat\bdy)$ is both 
a stationary point and a zero of $\bdg$.
The first component of the output $1\times 2$ cell array {\tt Z} is the 
computed zero
\[  \hat\bdx ~=~ (
\mbox{\tt \scriptsize 0.000000000000000},~
\mbox{\tt \scriptsize 0.000000000000000},~
\mbox{\tt \scriptsize 0.499435807628269},~
\mbox{\tt \scriptsize 2.002259318867864})
\]
that is accurate to the 16th digits with residual near the unit roundoff.
The Jacobian $\bdg_{\bdx\bdy}(\hat\bdx,\hat\bdy)$ is of nullity 1 that
is identical to the dimension of the zero $(\hat\bdx,\hat\bdy)$ to the
expanded mapping $\bdg$.
Namely, the ultrasingularity of the mapping $\bdf$ is deflated by expanding 
it to $\bdg$, resulting in a semiregular zero curve in $\R^4\times\R^4$ 
whose first component is the zero curve of $\bdf$.
\end{example}

At the current stage, theories of the depth-deflation approach 
(\ref{dfm}) for ultrasingular nonisolated zeros are still in development.
It is proved in \cite{HauWam13} that the Hauenstein-Wampler strong-deflation 
method terminates in finite many steps. 
It is highly desirable in practical numerical computation for the 
depth-deflation approach (\ref{dfm}) to terminate with the number of deflation 
steps bounded. 
The theories in \cite{DLZ,HauWam13} lead to a conjecture:
{\em Assume $\bdx_*$ is a $k$-dimensional ultrasingular zero of an analytic 
mapping $\bdf$ whose Jacobian $J(\bdx)$ maintains a constant nullity 
$n > k$ for all $\bdx\in\Omega\cap\bdf^{-1}(\bdo)$ where $\Omega$ 
is an open neighborhood of $\bdx_*$. 
Then the recursive depth-deflation process {\em (\ref{dfm})} terminates in 
finitely many steps so that $\bdx_*$ is a component of a semiregular 
$k$-dimensional zero of the final expanded system.}

\section{Conclusion}

Newton's iteration is not limited to solving for isolated regular solutions of 
nonlinear systems.
Nonisolated solutions of $\bdf(\bdx)\,=\,\bdo$ can be 
solved by the rank-$r$ Newton's iteration 
\[
\bdx_{k+1} ~~=~~ \bdx_k - J_\rkr(\bdx_k)^\dagger\,\bdf(\bdx_k) 
~~~~\mbox{for}~~~~ k\,=\,0,1,\ldots
\]
that locally quadratically converges to a semiregular zero of positive 
dimension where $J_\rkr(\bdx_k)^\dagger$ is the Moore-Penrose inverse of
the rank-$r$ approximation to the Jacobian $J(\bdx_k)$.
When the system is given with empirical data, the rank-$r$ Newton's iteration
still converges locally to a stationary point that accurately approximates
a nonisolated solution of the underlying exact system with an accuracy in the 
same order of the data, even if the solution of the given system is 
substantially altered by the data perturbation.

{\bf Acknowledgment.} The author thanks Dr. Tsung-Lin Lee for 
valuable suggestions and corrections on the manuscript.
The author also thanks Dr. Jonathan Hauenstein for pointing out the
deflation method and results in \cite{HauWam13}.

\providecommand{\bysame}{\leavevmode\hbox to3em{\hrulefill}\thinspace}
\providecommand{\MR}{\relax\ifhmode\unskip\space\fi MR }
\providecommand{\MRhref}[2]{%
  \href{http://www.ams.org/mathscinet-getitem?mr=#1}{#2}
}
\providecommand{\href}[2]{#2}

\end{document}